\documentclass[a4paper]{amsart}
\usepackage[utf8]{inputenc}

\usepackage{amsmath,amssymb}
\usepackage{lmodern,textcomp}
\usepackage{esint}  
\usepackage{enumerate} 
\usepackage{bbm}
\usepackage{mathrsfs}
\usepackage{upgreek}
\usepackage{cite}
\usepackage{graphicx} 
 
\usepackage{xcolor}
\definecolor{citation}{rgb}{0.2,0.55,0.2}
\definecolor{formula}{rgb}{0.1,0.2,0.6}
\definecolor{urlc}{rgb}{0,0,0.45}
\usepackage{hyperref}

\usepackage{a4wide}

\usepackage{amsfonts,amssymb,amsmath,amscd,amstext}
\usepackage{esint}
\usepackage{enumerate}
\usepackage{cancel,cite}
\usepackage{bbm,bm}

\newtheorem{theorem}{Theorem}[section]
\newtheorem{prop}[]{Proposition}[section]

\newtheorem{lemma}{Lemma}[section]

\newtheorem{defn}{Definition}[section]
\newtheorem{rem}{Remark}[section]

\numberwithin{equation}{section}

\DeclareMathAlphabet{\mathdutchcal}{U}{dutchcal}{m}{n}

\def\vs{\vspace{1mm}}

\def \eps{{\varepsilon}}
\def \r{{\mathbb{R}}}
\def \h{{\mathbb{H}^n}}

\def\mea{\mathcal{M}(\overline{\Omega})}

\def\XXint#1#2#3{{\setbox0=\hbox{$#1{#2#3}{\int}$}
		\vcenter{\hbox{$#2#3$}}\kern-.5\wd0}}

\def\ap{``}

\renewcommand{\epsilon}{\varepsilon}

\newcommand{\Ssub}{\Sob_\eps}
\newcommand{\Som}{\Sob_\Om}
\newcommand{\Sob}{S^\ast}

\newcommand{\rr}{\varrho}
\newcommand{\snr}[1]{\lvert #1\rvert}
\newcommand{\Dh}{\nabla_{X}}
\newcommand{\tows}{\stackrel{\ast}{\rightharpoonup}}

\newcommand{\Om}{\Omega}
\newcommand{\Omb}{\overline{\Omega}}

\begin{document}
	\title[Critical Sobolev embedding in the Heisenberg group]{Asymptotic approach to singular solutions\\ for the CR~Yamabe equation}
	
\author[G. Palatucci]{Giampiero Palatucci}  
\address[G. Palatucci]{Dipartimento di Scienze Matematiche, Fisiche e Informatiche
\newline\indent Universit\`a di Parma
\newline\indent Campus - Parco Area delle Scienze 53/A, 43124 Parma, Italy}  \email{\url{giampiero.palatucci@unipr.it}}

\author[M. Piccinini]{Mirco Piccinini} 
\address[M. Piccinini]{Dipartimento di Matematica
\newline\indent Universit\`a di Pisa
\newline\indent L.go~B.~Pontecorvo~5, 56127, Pisa, Italy
\newline\indent Dipartimento di Scienze Matematiche, Fisiche e Informatiche
\newline\indent Universit\`a di Parma
\newline\indent Campus - Parco Area delle Scienze 53/A, 43124 Parma, Italy}
\email{\url{mirco.piccinini@dm.unipi.it}}

\subjclass[2020]{35R03, 46E35, 35B33, 35J08, 35A15}

\keywords{Folland-Stein-Sobolev embeddings, Heisenberg group, CR Yamabe equation, concentration phenomena, Green functions}
	
\begin{abstract}
We study the loss of compactness for optimal functions of a subcritical approximations of the sharp Folland-Stein-Sobolev embedding in the Heisenberg group. When blow-up is prevent at the boundary we show that centered-symmetric maximizers in a Kor\'any ball concentrate at its center. Moreover, we obtain the exact blow-up rate and the pointwise asymptotic profile away from the pole, in turn establishing a natural counterpart of the classical Brezis and Peletier asymptotics for nearly critical Sobolev problems. 
\end{abstract}

\maketitle

\section{Introduction}
Let~$\h := (\mathbb{C}^n\times\mathbb{R},\circ,\{\delta_\lambda\}_{\lambda>0})$ be the usual Heisenberg-Weyl group, and let~$\mathring{\mathcal{D}}^{1,2}_X(\h)$ denote the completion of~$C^\infty_0(\h)$ with respect to the horizontal gradient norm~$\|\Dh \cdot\|_{L^2(\h)}$. In~\cite{FS74}, Folland and Stein proved that the following Sobolev-type inequality holds:
\begin{equation}\label{folland}
 \|u\|^{2^\ast}_{L^{2^\ast}(\h)} \leq \Sob \|\Dh u\|^{2^\ast}_{L^2(\h)} \qquad \forall u \in \mathring{\mathcal{D}}^{1,2}_X(\h),
\end{equation}
where~$\Sob$ is a positive constant and
$$
2^\ast:=\frac{2Q}{Q-2}
$$
is the Folland-Stein-Sobolev critical exponent, with~$Q:=2n+2$ denoting the homogeneous dimension associated with the anisotropic dilations~$\{\delta_\lambda\}_{\lambda>0}$; see, e.~\!g.,~\eqref{def_philambda} below.

The critical Sobolev inequality~\eqref{folland} has been an attractive object of study for the last decades, since it is inextricably linked to the lack of compactness of the related critical (Folland-Stein-)Sobolev embedding, and to the correspondent Euler-Lagrange equation in turn describing the important CR Yamabe problem~(\!\cite{JL87,Gam01}).

For what concern the strictly related study of the optimality of~\eqref{folland}, several results in accordance with the classical critical inequality in the Euclidean framework have been proven, despite the difficulties given by the sub-Riemannian geometry of the Heisenberg group~$\h$. On the contrary, several (somewhat expected) results are still open for the same reason; that is, the substantial difference with respect to the Euclidean framework in view of the complex non-commutative structure. 
The literature is too wide to attempt any comprehensive treatment in a single paper. We refer the interested readers to the very important papers~\cite{JL88,GL92,LU98,GV00,CU01,Loi05}, the recent book~\cite{IV11}, and the references therein.

\vspace{2mm}

In the present paper, we are interested into investigating some of the effects of the lack of compactness in the critical Sobolev embedding~\eqref{folland}, by analyzing the asymptotic behavior of the natural subcritical approximation of the Sobolev quotient.

\vspace{2mm}
Consider the following maximization problem,
\begin{equation}\label{critica0}
\Sob:=\sup\left\{\,\int_{\h}|u(\xi)|^{2^\ast}\,{\rm d}\xi \, : \, u\in \mathring{\mathcal{D}}^{1,2}_X(\h), \int_{\h}|\Dh u(\xi)|^2{\rm d}\xi \leq 1\right\}.
\end{equation}
The validity of~\eqref{folland} is equivalent to show that the constant~$\Sob$ defined in the display above is finite. The existence of the maximizers in~\eqref{critica0} is a difficult problem because of the intrinsic dilations and translations invariance of such inequality, as it analogously happens for the classical critical Sobolev inequality. The situation here is even more delicate because of the underlying non-Euclidean geometry of the Heisenberg group, and the  obstacles due to the related non-commutativity. The explicit form of the maximizers has been presented, amongst other results, in the breakthrough paper by Jerison and Lee~\cite{JL88}, together with the computation of the optimal constant in~\eqref{critica0}.
We also refer to the fundamental paper~\cite{FL12} where sharp constants for inequalities on~$\h$ have been derived for even more general cases, in turn obtaining sharp constants for the corresponding duals, which are the Sobolev inequalities for the sub-Laplacian and the conformal fractional Laplacians. 
\vspace{2mm}

For any bounded domain~$\Omega\subsetneq\h$, consider now the following Sobolev embedding in the same variational form as the one in~\eqref{critica0},
\begin{equation}\label{critica}
\Som:=\sup\left\{\,\int_{\Omega}|u(\xi)|^{2^\ast}\,{\rm d}\xi \, : \, u\in \mathring{\mathcal{D}}^{1,2}_X(\Omega), \int_{\Omega}|\Dh u(\xi)|^2{\rm d}\xi \leq 1\right\},
\end{equation}
where the Folland-Stein-Sobolev space~$\mathring{\mathcal{D}}^{1,2}_X(\Omega)$ is given as the closure of~$C^\infty_0(\Omega)$ with respect to the $L^2$-norm of the horizontal gradient in~$\Omega$.

\vspace{2mm}
One can check that~$\Som\equiv \Sob$ via a standard scaling argument on compactly supported smooth functions. For this, in view of the explicit form of the optimal functions in~\eqref{critica0} -- see forthcoming~Theorem~\ref{thm_optimal} -- the variational problem~\eqref{critica} has no maximizers. The situation changes considerably for the subcritical embeddings. Indeed, since~$\Omega$ is bounded,  the embedding~$\mathring{\mathcal{D}}^{1,2}_X(\Omega) \hookrightarrow L^{2^\ast-\eps}(\Omega)$ is compact (for any~$0 <\eps<2^\ast-2$), and this does guarantee the existence of a maximizer~$u_\eps\in \mathring{\mathcal{D}}^{1,2}_X(\Omega)$ for the related variational problem
   \begin{equation}\label{sobolev}
\Sob_\eps:=\sup\left\{\,\int_{\Omega}|u(\xi)|^{2^\ast-\eps}\,{\rm d}\xi \, : \, u\in \mathring{\mathcal{D}}^{1,2}_X(\Omega), \int_{\Omega}|\Dh u(\xi)|^2{\rm d}\xi \leq 1\right\}.
\end{equation}

Such a dichotomy is evident in the Euler-Lagrange equation for the energy functionals in~\eqref{sobolev}; that is,
\begin{equation}\label{equazione}
-\Delta_{X}   u_\eps = \lambda |u_\eps|^{2^\ast-\eps-2} u_\eps \, \ \text{in} \ \mathring{\mathcal{D}}^{1,2}_X(\Omega)',
\end{equation}
where~$\lambda$ is a Lagrange multiplier. Whereas when~$\eps>0$ the problem above has a solution~$u_\eps$, it becomes very delicate when~$\eps=0$: one falls in the aforementioned CR Yamabe equation, and even the existence of solutions is not granted. In particular, the existence and various properties of the solutions do strongly depend on the geometry and the topology of the domain~$\Omega$. We refer for instance to:~\cite{LU98,Ugu99}
for nonexistence of nonnegative solutions when~$\Om$ is a certain half-space; \cite{CU01}~where the authors show the existence of a solution in the case when the domain~$\Omega$ has at least a nontrivial suitable homology group; 
 \cite{GL92} for existence and nonexistence results for even more general nonlinearity.
\vspace{2mm}

In view of such a qualitative change when~$\eps=0$ in both~\eqref{sobolev} and~\eqref{equazione}, it seems natural to analyze the asymptotic behaviour as~$\eps$ goes to~$0$ of  the corresponding optimal functions~$u_\eps$ of the embedding $\mathring{\mathcal{D}}^{1,2}_X(\Omega) \hookrightarrow L^{2^\ast\!-\eps}(\Omega)$. This is the aim of the present paper.
\vspace{2mm}

For what concerns the~Euclidean~counterpart of such an investigation, several results have been obtained, mostly via fine estimates and a standard regularity elliptic approach of the special class of solutions of the equation~\eqref{equazione} being maximizers for the related Sobolev embedding.

 \vspace{2mm}
 On the contrary, for what concerns the~Heisenberg~panorama the scene is basically empty  in view of the many difficulties  naturally arising in such a framework. Indeed, the non-commutative group structure precludes the free generalization of several tools such as
 symmetric decreasing rearrangements, ODEs techniques as well as regularity approximations. 
Nevertheless, some fundamental results available in the Euclidean framework where successfully extended in~$\h$: the Bahri and Coron conjecture in~\cite{BC88} in the aforementioned paper~\cite{CU01}, as well as the  non-existence criteria on several relevant class of proper subset of~$\h$ in~\cite{LU98,GL92}.
\vspace{2mm}

Very recently, it has been proven in~\cite{PPT25} that, up to subsequences, optimal functions~$u_\eps$ for the subcritical Sobolev embedding~\eqref{sobolev} do concentrate horizontal energy at one point~$\xi_{\rm o} \in \Omb$, even without  prescribing any regularity assumptions nor special geometric features on the domain~$\Omega$, in clear accordance with their Euclidean counterpart in~\cite{FM99,AG03,Pal11,Pal11b}; see also the recent results in the nonlocal settings in~\cite{PT23}.

\begin{theorem}[Theorem~1.2 in~\cite{PPT25}]\label{cor_concentration}
Let~$\Om\subsetneq \h$ be a bounded domain, let~$\mea$ denote the family of positive Radon measures on~$\Omb$, and let~$U_\eps\in \mathring{\mathcal{D}}^{1,2}_X(\Omega)$ be a maximizer for~$\Sob_\eps$, normalized by
$$
\int_{\Omega}|\Dh U_\eps(\xi)|^2\,{\rm d}\xi=1.
$$
Then, for every sequence~$\eps_k\to0^+$, up to subsequences, there exists~$\xi_{\rm o}\in\Omb$ such that
$$
U_{\eps_k}\rightharpoonup0
\quad\text{in }L^{2^\ast}(\Omega),
$$
and
$$
|\Dh U_{\eps_k}|^2\,{\rm d}\xi
\tows
\boldsymbol{\delta}_{\xi_{\rm o}}
\quad\text{in }\mea.
$$
\end{theorem}

\vspace{2mm}

Now, a natural question arises: can we  localize the blow up; i.~\!e., is the concentration point~$\xi_{\rm o}$ related to some extent to the geometry of the domain~$\Omega$\,?
 \\ For classical elliptic equations with critical growth nonlinearities,  Atkinson and Peletier, via ODEs methods, proved in~\cite{AP87} that the blow up (as~$\eps \searrow 0$) of solutions to the critical equation
 	\begin{equation}\label{eq:el-sob}
	\begin{cases}
 		-\Delta u = u^{\frac{n+2}{n-2}-\eps} & \text{in}~\Om\\
		u > 0 & \text{in}~\Om\,,\\
		u =0 & \text{in}~\partial\Om\,,
	\end{cases}
 	\end{equation}
 	in the case when~$\Om \subset \r^3$ does coincide with the unit ball, satisfy
 	\[
 	\lim_{\eps \to 0^+} \eps u_\eps^2(0) = \frac{32}{\pi} \quad \text{and} \quad  \lim_{\eps \to 0^+} \frac{ u_\eps(x)}{\sqrt{\eps}} =\frac{1}{4}\sqrt{\frac{\pi}{2}}\left(\frac{1}{\snr{x}}-1\right).
 	\]
 	Subsequently, by relying on purely variational techniques, Brezis and Peletier  in~\cite{BP89}  extended such results to the case of spherical domains. In particular, 
	 it has been showed that  subcritical solutions concentrate at one special point of the domain. Moreover, the authors also conjectured that an analogous result should hold for non spherical	domains and for higher dimensions as well. This conjecture was later proved to be true, independently, in the case of smooth domains by Han~\cite{Han91} and Rey~\cite{Rey89} by showing that the solutions of~\eqref{eq:el-sob}, with maximal Sobolev energy, concentrate energy at one point which can be localized via the Green's function associated with the underlying domain.
 
\vspace{2mm}
  The involved proofs strongly rely {as well as on various available techniques in the Euclidean framework such as,~e.~\!g., moving planes method, Kelvin transform, etc..., also} on {the availability  of various boundary} regularity {results for standard elliptic equations}. {This last feature} is in clear contrast with the complexity {faced in the present paper}. As well known, even if the domain~$\Omega$ is smooth, the situation is dramatically different because of the possible presence of characteristic points on the boundary~$\partial\Omega$. 
  At such points the vector fields forming the  principal part of the relevant operator~$\Delta_{X}  $ become tangent to the boundary. Hence, near those characteristic points -- as firstly discovered by Jerison~\cite{Jer81} -- even harmonic functionscan encounter a sudden loss of regularity. {Indeed, Jerison built an explicit solution in the domain $\big\{\xi=(z,t) \in \h: \, t>-M\snr{z}^2\big\}$, for a suitably choice fo $M>0$, vanishing on the boundary and {having at most H\"older regularity} near its isolated characteristic point~${0}$.}
   Also, one did not want to work in the restricted class of domains not having characteristic points; that is, by still including interesting sets as e.~\!g. the torus obtained by revolting the sphere~$\mathbb{S}^{2n}$ around the $t$-axis~\cite{AG21}, but unfortunately excluding an extremely  wide class of regular domains  which play a pervasive role in several relevant problems in the Heisenberg group, as e.~\!g. the level sets of the Jerison and Lee extremal functions~\eqref{talentiane_2} and those of the Folland fundamental solution; i.~\!e., the Kor\'anyi balls. {Nevertheless, some important results have been obtained for maximizing sequence of~$\Ssub$ in non-characteristic domain. Indeed, in~\cite{MMP13}, the authors are able to construct a concentrating sequence of solutions for certain non-degenerate critical point of  the regular part of the Green function of~$\Om$. We also refer  to the references in the aforementioned paper.}

\vspace{2mm}
The purpose of the present paper is to go further Theorem~\ref{cor_concentration} and to determine the quantitative asymptotics of its concentration.
The canonical model is the Kor\'anyi ball
$$
B_R(0):=\big\{\xi=(x,t)\in\h:\ \big(|x|^4+t^2\big)^\frac{1}{4}<R\big\}.
$$
This is the natural counterpart of the Euclidean ball: it is defined by the homogeneous norm associated with the fundamental solution~$|\xi|_{\mathbb{H}}:=\big(|x|^4+t^2\big)^\frac{1}{4}$, it is invariant under the unitary symmetries in the horizontal variables and under the reflection~$(z,t)\mapsto(\overline z,-t)$. In this setting, both the Green and Robin function are explicit
$$
G_{B_R(0)}(\xi;0) :=
\frac{1}{C_Q}\left(\frac{1}{|\xi|_{\mathbb H}^{Q-2}}-\frac{1}{R^{Q-2}}\right) =: K(\xi)-\mathcal{R}_{B_R(0)}(0)\,,
$$
where~$C_Q$ is the dimensional constant in the Folland fundamental solution~$K(\xi)=C_Q^{-1}|\xi|_{\mathbb H}^{2-Q}$.

\vspace{2mm}
In order to state our main results, we introduce some further notation.

\begin{defn}
Let~$U_\eps$ be a nonnegative maximizer in~\eqref{sobolev}, so that~$\|\Dh U_\eps\|_{L^2(\Om)}=1$.
Then, we call scalar normalized maximiser the function
\begin{equation}\label{eq:scalar-change-intro}
u_\eps
:=
(\Sob_\eps)^{-\frac{1}{2^\ast-\eps-2}}U_\eps.
\end{equation}
\end{defn}

From~\eqref{eq:scalar-change-intro}, it follows that~$u_\eps$ solves the Dirichlet problem
\begin{equation}
\begin{cases}\label{eq:normalized-critical-subcritical}
-\Delta_Xu_\eps=u_\eps^{\frac{Q+2}{Q-2}-\eps}
& \text{in}~\Om\,,\\
u_\eps=0 & \text{on}~\partial\Omega.
\end{cases}
\end{equation}
Moreover,
\begin{align}
& (\Sob_\eps)^{-\frac{2}{2^\ast-\eps-2}}
=
\int_\Omega |\Dh u_\eps|^2\,{\rm d}\xi
=
\int_\Omega u_\eps^{2^\ast-\eps}\,{\rm d}\xi .
\label{eq:scalar-normalized-mass-intro}
\end{align}

We shall use that
\begin{equation}\label{eq:Sob-eps-limit}
\Sob_\eps \to \Sob
\qquad\text{as }\eps\to0^+\,,
\end{equation}
as can be proved by a simple application of H\"older's Inequality; see, e.~\!g.,~\cite[Proposition~2.3]{PPT25}. In particular,
\begin{equation}\label{eq:def-Sstar-normalized}
(\Sob)^{-\frac{Q-2}{2}}
=
\lim_{\eps\to0^+}
\int_\Omega u_\eps^{2^\ast-\eps}\,{\rm d}\xi
=
\lim_{\eps\to0^+}
\int_\Omega |\Dh u_\eps|^2\,{\rm d}\xi.
\end{equation}

For the scalar-normalized sequence in~\eqref{eq:scalar-change-intro}, the concentration statement reads, along the selected subsequence,
\begin{equation}\label{eq:scalar-normalized-concentration-measure}
|\Dh u_\eps|^2\,{\rm d}\xi
\tows
(\Sob)^{-\frac{Q-2}{2}}\,\boldsymbol{\delta}_{\xi_{\rm o}}
\qquad
\text{in }\mea.
\end{equation}
Let~$U$ denote the normalized Jerison--Lee extremal satisfying
$$
-\Delta_XU=U^{2^\ast-1}
\quad\text{in }\h,
\qquad
U(0)=\max_{\h}U=1,
$$
and set
$$
\bar c_n:=\int_\h U^{2^\ast-1}\,{\rm d}\xi.
$$

The main result of the paper is the following explicit Atkinsons--Brezis--Peletier type asymptotic formula in the centered Kor\'anyi ball.

\begin{theorem}\label{cor:centered-koranyi-ball}
Let~$u_\eps$ be a family of scalar-normalized maximizers in~$B_R(0)$. Assume that~$u_\eps$ is centered-symmetric, i.~\!e.
$$
u_\eps(Az,t)=u_\eps(z,t)
\quad\forall A\in\mathcal U(n),
\qquad
u_\eps(\overline z,-t)=u_\eps(z,t),
$$
and that there exist~$\delta\in(0,R)$,~$C>0$, and~$\eps_{\rm o}>0$ such that
\begin{equation}\label{eq:boundary-exlusion}
\sup_{\xi\in B_R(0)\setminus B_{R-\delta}(0)}u_\eps(\xi)
\leq C
\qquad
\forall\eps\in(0,\eps_{\rm o}).
\end{equation}
Then the following assertions hold true:
\begin{enumerate}[{\it (i)}]
\item $u_\eps$ concentrates at the center of the ball, i.~\!e.
$$
|\Dh u_\eps|^2\,{\rm d}\xi
\tows
(\Sob)^{-\frac{Q-2}{2}}\,\boldsymbol\delta_0.
$$
\item One has
\begin{equation}\label{eq:blowup-koranyi-ball}
\lim_{\eps\to0^+}
\eps\|u_\eps\|_{L^\infty(B_R(0))}^2
=
c_1(n)\mathcal R_{B_R(0)}(0),
\end{equation}
where
$$
c_1(n)
:=
\frac{2^\ast\bar c_n^2}{(\Sob)^{-\frac{Q-2}{2}}}.
$$
\item For every~$\xi\in B_R(0)\setminus\{0\}$,
\begin{equation}\label{eq:profile-koranyi-ball}
\lim_{\eps\to0^+}\frac{u_\eps(\xi)}{\sqrt\eps}
=
c_2(n)
\left(
\frac{R^{\frac{Q-2}{2}}}{|\xi|_{\mathbb H}^{Q-2}}
-
\frac{1}{R^{\frac{Q-2}{2}}}
\right),
\end{equation}
where
$$
c_2(n)
:=
\sqrt{\frac{(\Sob)^{-\frac{Q-2}{2}}}{2^\ast C_Q}}.
$$
\end{enumerate}
\end{theorem}

\begin{rem}\label{rem:normalization}
{\rm
The constant in~\eqref{eq:blowup-koranyi-ball} depends on the scalar normalization chosen for the maximizing sequence. If one works instead with the variational normalization
$$
\int_\Omega |\Dh U_\eps|^2\,{\rm d}\xi=1,
$$
then the same statement holds after applying the scalar change from $U_\eps$ to the solution $u_\eps$ of~\eqref{eq:normalized-critical-subcritical}.}
\end{rem}

\begin{rem}\label{rem:other-approximations}
{\rm
As in~\cite{Han91}, the argument applies more generally to positive solutions of~\eqref{eq:normalized-critical-subcritical} whose Sobolev quotient is asymptotically sharp. One may also consider other subcritical perturbations of the critical CR Yamabe equation, for instance families of the form
$$
-\Delta_Xu_\eps
=
\lambda u_\eps^{2^\ast-1}
+\eps u_\eps,
$$
where the loss of compactness again appears as $\eps\to0$. We do not pursue this direction here, in order to focus on the optimal functions for the subcritical Folland-Stein-Sobolev embeddings.
}
\end{rem}

\vspace{2mm}

A key ingredient in the proof of Theorem~\ref{cor:centered-koranyi-ball} is a sharp asymptotic control of the scalar-normalized maximizing sequence by the Jerison-Lee extremals, which is our second main result. This estimate is global and is independent of the explicit Green-function computation in the ball.

\begin{theorem}\label{han}
Let~$\Omega\subsetneq\h$ be a bounded smooth domain. Let $u_\eps\in \mathring{\mathcal{D}}^{1,2}_X(\Omega)$ be the scalar-normalized maximizer solving~\eqref{eq:normalized-critical-subcritical}. Assume there exist~$\delta>0$,~$C_0>0$, and~$\eps_{\rm o}>0$ such that
$$
\sup_{\{\xi\in\Om:\ \operatorname{dist}_{\mathbb R^{2n+1}}(\xi,\partial\Om)<\delta\}}u_\eps(\xi)
\leq C_0
\qquad
\forall\eps\in(0,\eps_{\rm o}).
$$
Set
$$
M_\eps:=\|u_\eps\|_{L^\infty(\Omega)}, \quad u_\eps(\eta_\eps)=M_\eps \quad \text{and} \quad \rho_\eps := M_\eps^{-\frac{2^\ast-\eps-2}{2}}.
$$
Then, when~$\eps\to0$, up to subsequences, there exists a point~$\xi_{\rm o}\in\Omega$ such that
$$
|\Dh u_\eps|^2\,{\rm d}\xi
\tows
(\Sob)^{-\frac{Q-2}{2}}\,\boldsymbol{\delta}_{\xi_{\rm o}}
\quad\text{in }\mea,
$$
and the following assertions hold true:
\begin{enumerate}[(i)]
\item $M_\eps \to +\infty$;
\item $\eta_\eps\to \xi_{\rm o}$;
\item Letting~$U$ be the normalized Jerison-Lee extremal
$$
\begin{cases}
-\Delta_XU=U^{2^\ast-1}
\quad \text{in}~\mathbb H^n,\\
U(0)=\max\limits_{\mathbb H^n}U=1,
\end{cases}
$$
and defining
$$
v_\eps(\zeta):=M_\eps^{-1} u_\eps\big(\eta_\eps\circ\delta_{\rho_\eps}(\zeta)\big),
\qquad
\zeta\in\Omega_{\eta_\eps,\rho_\eps},
$$
where
$$
\Omega_{\eta_\eps,\rho_\eps}:=
\delta_{\rho_\eps^{-1}}\big(\eta_\eps^{-1}\circ\Omega\big),
$$
it holds that~$v_\eps \to U$ locally uniformly in~$\mathbb H^n$ and locally in~$\mathring{\mathcal{D}}^{1,2}_X(\h)$.
\end{enumerate}
Furthermore, there exists a constant
$C=C(n,\Omega,\delta,C_0)>0$, independent of~$\eps$, such that
\begin{equation}\label{bound_max_seq}
u_\eps(\xi)
\leq
C\,M_\eps\,
U\left(
\delta_\frac{1}{\rho_\eps}
\big(\eta_\eps^{-1}\circ\xi\big)
\right)
\quad
\text{for every}~\xi\in\Omega .
\end{equation}
\end{theorem}

The estimate above is the Heisenberg analogue of the global bubble control used by Han in the Euclidean proof~\cite{Han91}. Here the argument has to account for the non-commutative translations and for the fact that the Jerison-Lee extremals cannot be reduced, after arbitrary translations and dilations, to functions of the Kor\'anyi gauge alone. The proof combines the concentration theorem recalled above, the intrinsic scaling of the sub-Laplacian, a sharp-quotient argument based on the Brezis--Lieb lemma, the $\mathbb H$-Kelvin transform, and boundary estimates adapted to the characteristic geometry of the Kor\'anyi ball.

\subsection{Related open problems and further developments}

Starting from the results proved in the present paper, several questions naturally arise.

\begin{itemize}

\item One may investigate fractional counterparts of the present asymptotics, replacing the first-order Folland-Stein norm in~\eqref{folland} by a fractional norm of order $s\in(0,1)$ and the critical exponent by $2_s^\ast=2Q/(Q-2s)$. In this direction we mention~\cite{GLV23}, the nonlinear tail estimates in~\cite{MPPP23}, related results in~\cite{PP22}, and the extension/scattering approach to conformally invariant fractional powers of the sub-Laplacian in~\cite{FGMT15}.\smallskip

\item It would be natural to ask to what extent the Green-profile asymptotics proved here in the centered Kor\'anyi ball persist under different geometric assumptions on the domain, still allowing characteristic boundary points. A possible framework is that of non-tangentially accessible domains satisfying an intrinsic outer ball condition, as introduced in~\cite{CGN02} for the Dirichlet problem with summable boundary data.\smallskip

\item The behavior of concentration points in non-smooth domains remains a delicate issue. In the Euclidean framework, Flucher, Garroni and M\"uller~\cite{FGM02} constructed a non-smooth domain whose Robin function attains its infimum at the boundary, and Pistoia and Rey~\cite{PR03} showed that concentration can occur at the boundary in such a domain. Understanding whether analogous phenomena can occur in the Heisenberg group would require controlling the interaction between boundary singularities, characteristic points, and the Green function; see also~\cite{Wei98}.\smallskip

\item Other possible developments include critical energies beyond the pure Sobolev quotient, problems involving the second critical exponent as in~\cite{Pas93,DPMP10}, and analogues of the third Brezis-Peletier conjecture~\cite{BP89}. The latter was recently solved in the Euclidean setting in~\cite{FKK23}, but its Heisenberg counterpart appears to require new tools because both the critical geometry and the boundary regularity theory are substantially more involved.
\end{itemize}\smallskip

The estimates and techniques developed here are intended to provide a basis for these further questions on the loss of compactness in critical Folland-Stein-Sobolev embeddings and in related CR Yamabe-type equations.

\subsection{Organization of the paper}

In Section~\ref{sec_preliminaries} we fix the notation and recall the analytic and geometric tools used throughout the paper: the structure of the Heisenberg group, the Folland fundamental solution, characteristic boundary points, intrinsic H\"older spaces, the Pohozaev identity, and the Jerison-Lee extremals. In Section~\ref{sec_asymptotic} we prove the sharp asymptotic control of subcritical maximizers by Jerison-Lee bubbles, culminating in Theorem~\ref{han}. In Section~\ref{sec_localization} we prove the centered Kor\'anyi-ball asymptotics of Theorem~\ref{cor:centered-koranyi-ball}, combining the Green-function identity, the boundary estimates near the characteristic set, and the anisotropic Pohozaev identity.

\subsection*{Acknowledgements}
The authors are supported by INdAM projects \ap Problemi non locali: teoria cinetica e non uniforme ellitticit\`a'' -- CUP\_E53C22001930001, and
\ap Local vs Nonlocal: mixed type operators and nonuniform ellipticity'' -- CUP\_D91B21005370003. 
The first author is also supported by PRIN 2022 Project ``Geometric Evolution Problems and Shape Optimization (GEPSO)'' -- CUP\_D53D23005820006, and by PRIN 
2022 PNRR ``Magnetic skyrmions, skyrmionic bubbles and domain walls for spintronic applications'' -- CUP\_D53D23018980001, 
PNRR Italia Domani, financed by EU via NextGenerationEU.

The results in this paper have been announced in the preliminary research report~\cite{PP23}.

\medskip

\section{Preliminaries}\label{sec_preliminaries}

In this section we fix the notation and recall the basic tools used throughout the paper
{ such as the
the 
structure of the Heisenberg group, the relevant function spaces and 
some boundary Schauder-type estimates and integral identities.}

\subsection{The Heisenberg group}
Throughout the paper we denote points in~$\r^{2n+1}$ by~$\xi:=(x,t)=(x_1,\dots,x_{2n},t)\in\r^{2n}\times\r
$.
Equivalently, we often write~$z:=(z_1,\dots,z_n)$ with~$z_j:=x_j+i x_{n+j}$, so that~$\h$ is identified with~$\mathbb C^n\times\r$. The Heisenberg group~$\h$ is the analytic, simply connected $(2n+1)$-dimensional Lie group whose Lie algebra~$\mathfrak g$ admits the stratification
\[
\mathfrak{g}=V^1\oplus V^2,
\quad
[V^1,V^1]=V^2,
\quad
[V^1,V^2]=\{0\}.
\]
A basis of left-invariant vector fields spanning the first layer~$V^1$ is given by
\[
X_j:=\partial_{x_j}+2x_{n+j}\partial_t,
\qquad
X_{n+j}:=\partial_{x_{n+j}}-2x_j\partial_t,
\qquad 1\leq j\leq n.
\]
Moreover, since the only non-trivial commutators are
\[
[X_j,X_{n+j}]:=X_jX_{n+j}-X_{n+j}X_j=-4\partial_t,
\qquad 1\leq j\leq n,
\]
a basis for the second layer~$V^2$ is 
given
by~$
T:=\partial_t$, so that the stratification of the Lie algebra~$\mathfrak{g}$ becomes
\[
\mathfrak{g}=
\operatorname{span}\{X_1,\dots,X_{2n}\}
\oplus
\operatorname{span}\{T\}.
\]

The Heisenberg group~$\h$ can be identified with the triple~$(\r^{2n+1},\circ,\{\delta_\lambda\}_{\lambda>0})$, where~$\circ$ is the polynomial group multiplication law 
\[
({x},t)\circ({y},s) :=
\bigg(
x_1+y_1,\dots,x_{2n}+y_{2n},
t+s+2\sum_{j=1}^n(x_{n+j}y_j-x_jy_{n+j})
\bigg),
\]
for every~$(x,t),(y,s)\in\h$. With this convention the inverse of a point~$\xi=(x,t)$ is~$\xi^{-1}=(-x,-t)$ and the group identity is given by~$0:=(0,0)$.

The group of anisotropic dilations~$\{\delta_\lambda\}_{\lambda>0}$ is 
\begin{equation}\label{def_philambda}
\xi \mapsto\delta_\lambda(\xi):=(\lambda x,\lambda^2t),
\quad \lambda>0.
\end{equation}
Given~$\xi'\in\h$, the left translation by~$\xi'$ is denoted by
\begin{equation}\label{def_tau}
\tau_{\xi'}(\xi):=\xi'\circ\xi,
\qquad \xi\in\h .
\end{equation}
As  customary, we denote with~$Q:=\dim V^1+2\dim V^2=2n+2$ the homogeneous dimension associated with the dilations~\eqref{def_philambda}.

We denote by~$
\Dh:=(X_1,\dots,X_{2n})$ 
the horizontal (or intrinsic) gradient~in~$\h$, while we indicate by~$
\nabla:=(\partial_{x_1},\dots,\partial_{x_{2n}},\partial_t)$ the standard Euclidean gradient. The Kohn sub-Laplacian~$\Delta_X$ is the second order operator left-invariant with respect to the translations in~\eqref{def_tau} and homogeneous of degree~$2$ with respect to the dilations in~\eqref{def_philambda}
$$
\Delta_X:=\sum_{j=1}^{2n}X_j^2 .
$$

It is well known, see, e.~\!g.~\cite{Fol75}, that~$\Delta_X$ admits a fundamental solution~$K\in C^\infty(\h\setminus\{0\})\cap L^1_{\rm loc}(\h)$,~$K(\xi) \to 0$ when~$\xi$ tends to infinity, such that~$\Delta_XK=-\boldsymbol{\delta}_0$ in the distributional sense, i.~\!e.
\[
\int_{\h}K(\xi)\Delta_X\phi(\xi)\,{\rm d}\xi
=
-\phi(0)
\quad
\forall \phi\in C^\infty_0(\h).
\]
\begin{defn}
    We call Gauge norm on~$\h$ a homogeneous symmetric norm~$d$ smooth out of the origin and such that
\[
\Delta_X(d(\xi)^{2-Q})=0 \quad \forall \xi \neq 0.
\]
\end{defn}

We use the standard gauge on~$\h$ also known as Kor\'anyi gauge
\[
|\xi|_{\mathbb{H}}:=\big(|x|^4+t^2\big)^{1/4},
\quad \forall\xi=(x,t)\in\h.
\]
In this way, we have that
\begin{equation}\label{eq:fundamental}
K(\xi)=C_Q^{-1}|\xi|_{\mathbb{H}}^{2-Q}
=
\frac{1}{C_Q\big(|x|^4+t^2\big)^{(Q-2)/4}},
\qquad \xi\neq0,
\end{equation}
where~$C_Q>0$ is an explicit dimensional constant. We define the fundamental solution with pole at~$\eta$ as~$K(\eta^{-1}\circ\xi)=C_Q^{-1}|\eta^{-1}\circ\xi|_{\mathbb{H}}^{2-Q}$.

As customary, we will denote with~$B_r(\eta)$ the gauge ball centered at~$\eta\in\h$ and with radius~$r>0$
given by~$
B_r(\eta):=\{\xi\in\h:\ |\eta^{-1}\circ\xi|_{\mathbb{H}}<r\}$.

\subsection{Smooth domains, characteristic points, and intrinsic H\"older spaces}

We say that a domain~$\Omega\subsetneq\h$ is of class~$C^\infty$ if for every~$\xi\in\partial\Omega$ there exist a neighborhood~$U_\xi$ of $\xi$ and a function~$\Upphi_\xi\in C^\infty(U_\xi)$ such that
\begin{eqnarray*}
U_\xi\cap\Omega & = & \{\eta\in U_\xi:\Upphi_\xi(\eta)<0\}\,,\\
U_\xi\cap\partial\Omega & = & \{\eta\in U_\xi:\Upphi_\xi(\eta)=0\},
\end{eqnarray*}
and
\[
\nabla\Upphi_\xi\neq 0
\quad \text{on}~U_\xi\cap\partial\Omega.
\]
A boundary point~$\xi\in\partial\Omega$ is called {\it characteristic} if~$\Dh\Upphi_\xi(\xi)=0$
and we denote with~$\Sigma(\Om)$ the set of all characteristic points of~$\Om$, i.~\!e.
\[
\Sigma(\Omega):=
\{\xi\in\partial\Omega:\Dh\Upphi_\xi(\xi)=0\}.
\]

We recall the definition of intrinsic H\"older spaces. 

\begin{defn}
  Let~$\beta\in(0,1)$ and~$\Omega\subseteq\h$. A function~$u:\Omega\to\r$ belongs to~$C^{0,\beta}(\Omega)$ if
$$
[u]_{C^{0,\beta}(\Omega)}
:= \sup_{\substack{\xi,\eta\in\Omega\\ \xi\neq\eta}}
\frac{|u(\xi)-u(\eta)|}
{|\eta^{-1}\circ\xi|_{\mathbb H}^{\beta}} <\infty .
$$
\end{defn}

The space~$C^{0,\beta}(\Om)$ is a Banach space once endowed with the norm
$$
\|u\|_{C^{0,\beta}(\Omega)}
:= \|u\|_{L^\infty(\Omega)} + [u]_{C^{0,\beta}(\Omega)}.
$$
{
For~$k\in\mathbb N$, we use the standard non-isotropic convention: if~$
Z=X_{i_1}\cdots X_{i_\ell}$
is an ordered differential monomial in the vector fields~$X_1,\dots, X_{2n}$, we assign to~$Z$ the homogeneous order
\[
{\rm ord}(Z):=\ell.
\]
We say that~$u\in C^{k,\beta}(\Omega)$ if all derivatives~$Zu$ with~${\rm ord}(Z)\leq k$ exist and are continuous in~$\Omega$, and all derivatives~$Zu$ with~${\rm ord}(Z)=k$ belong to~$C^{0,\beta}(\Omega)$.
}

We shall also use intrinsic H\"older spaces on non-characteristic portions of the boundary. Let~$\xi,\eta\in\partial\Omega$  we define the induced distance~$\widehat{d}$ on~$\partial\Om$ as
\[
\widehat{d}(\xi,\eta):=|\eta^{-1}\circ\xi|_{\mathbb{H}}.
\]
For~$r>0$, we write
\[
\widehat B_r(\xi):=B_r(\xi)\cap\partial\Omega.
\]
For~$\beta \in (0,1)$, we say that a continuous function~$u:\partial\Omega\to\r$ belongs to~$C^{0,\beta}(\partial\Omega)$ if
\[
[u]_{C^{0,\beta}(\partial\Omega)}
:=
\sup_{\substack{\xi,\eta\in\partial\Omega\\ \xi\neq\eta}}
\frac{|u(\xi)-u(\eta)|}{\widehat d(\xi,\eta)^\beta}
<\infty\,,
\]
and we endow~$C^{0,\beta}(\partial\Om)$ with the following norm,
$$
\|u\|_{C^{0,\beta}(\partial\Om)}:= \sup_{\xi \in \partial\Om}|u(\xi)|+[u]_{C^{0,\beta}(\partial\Om)}.
$$

We now recall the boundary spaces~$\Gamma^{k,\beta}$ introduced by Jerison in~\cite{Jer81}. Let~$\xi\in\partial\Omega$ be a non-characteristic point. Then the horizontal outward unit normal
\[
{n}^{\mathbb{H}}_{\partial\Om}
:=
\frac{\Dh\Upphi_\xi}{|\Dh\Upphi_\xi|},
\]
is well defined near~$\xi$. Therefore, on a sufficiently small non-characteristic neighborhood~$U_\xi$ of $\xi$, we can describe~$\partial\Om \cap U_\xi$ by intrinsic coordinates~$
\hat{\zeta}=(\zeta_1,\dots,\zeta_{2n}) \in \partial\Om \cap U_\xi$
obtained from the inverse of the exponential map associated with the adapted boundary frame; see~\cite[Section~3]{Jer81}. For a multi-index~$J=(j_1,\dots,j_{2n})$, set
\[
\deg(J):=j_1+\cdots+j_{2n-1}+2j_{2n},
\]
and
\[
\hat{\zeta}^J
:=
\zeta_1^{j_1}\cdots\zeta_{2n}^{j_{2n}}.
\]
A polynomial of intrinsic degree at most $k$ in the boundary coordinates is a polynomial of the form
\[
P(\hat{\zeta})
=
\sum_{\deg(J)\leq k}a_J\hat{\zeta}^J\,,
\]
for some constant~$a_J$.
If~$\eta$ is represented by the coordinates~$\hat{\zeta}$, we write~$P_\xi(\eta):=P(\hat{\zeta})$.

\begin{defn}\label{def:intr-hld}
Let~$\beta\in(0,1)$,~$k\in\mathbb N\cup\{0\}$ and let~$B_r(\eta_0)\subset\h$ be a ball such that~$\partial\Om\cap B_r(\eta_0)$ is non-characteristic. A bounded function~$u$ belongs to~$\Gamma^{k,\beta}(\partial\Om\cap B_r(\eta_0))$ if there exists~$C>0$ such that, for every~$\xi\in\partial\Om\cap B_r(\eta_0)$, there exists a polynomial~$P_\xi$ of intrinsic degree at most~$k$ satisfying
\begin{equation}\label{eq:intrinsic-hld-space}
|u(\eta)-P_\xi(\eta)|
\leq
C\widehat d(\eta,\xi)^{k+\beta}
\end{equation}
for every~$\eta\in\partial\Omega\cap B_r(\eta_0)$. Then
$$
\|u\|_{\Gamma^{k,\beta}(\partial\Om\cap B_r(\eta_0))}
:=
\sup_{\xi\in\partial\Om\cap B_r(\eta_0)}|u(\xi)|
+
[u]_{\Gamma^{k,\beta}(\partial\Om\cap B_r(\eta_0))},
$$
where~$[u]_{\Gamma^{k,\beta}}$ is the least admissible constant~$C$, together with the supremum of the coefficients of the polynomials~$P_\xi$ appearing in~\eqref{eq:intrinsic-hld-space}.
\end{defn}

{These classes are the non-isotropic H\"older classes (see, e.~\!g.,~\cite{FS82}) introduced by Jerison in~\cite[Section~4]{Jer81}.}

\begin{rem}{\rm
The same notation on~$\overline{\Omega \cap B_r(\eta_0)}$ will refer to the corresponding boundary non-isotropic H\"older class combined with the usual interior Folland-Stein H\"older control away from the boundary.}
\end{rem}

With this notation we recall the boundary Schauder-type estimate that will be used below. For further Schauder estimates at the boundary in Lie groups we refer to~\cite{CGS21,BGM22}.

\begin{theorem}[Theorem~1.1 in~\cite{BGM19}]\label{boundary_hld}
Let~$\Omega\subset\h$ be a bounded domain of class~$C^{1,\beta}$ for some~$\beta\in(0,1)$, and let~$B_r(\eta_0)$ be such that~$\partial\Omega \cap B_r(\eta_0)$ is contained in a non-characteristic portion of~$\partial\Omega$. Let~$
u\in H_X^1(\Omega\cap B_r(\eta_0))\cap C(\overline{\Omega\cap B_r(\eta_0)})$
be a weak solution of
\begin{equation}\label{eq:bgm}
\begin{cases}
-\Delta_Xu=f & \text{in } \Omega \cap B_r(\eta_0),\\
u=0 & \text{on}~\partial\Omega \cap B_r(\eta_0),
\end{cases}
\end{equation}
with~$f\in L^\infty(\Omega \cap B_r(\eta_0))$. 

Then,~$u\in\Gamma^{1,\beta}(\overline{\Omega \cap B_\frac{r}{2}(\eta_0)})$, 
and we have the apriori estimate
\[
\|\Dh u\|_{L^\infty(\Omega\cap B_\frac{r}{2}(\eta_0))}
+
r^\beta[\Dh u]_{C^{0,\beta}(\Om \cap B_\frac{r}{2}(\eta_0))}
\leq
\frac{C}{r}\big(
\|u\|_{L^\infty(\Omega \cap B_r(\eta_0))}
\ + r^2\|f\|_{L^\infty(\Omega\cap B_r(\eta_0))}
\big),
\]
where~$C$ depends only on~$n$,~$\beta$ and on~$\partial\Omega$.
\end{theorem}

\begin{rem}\label{rem:bgm-boundary}
{\rm
In particular, Theorem~\ref{boundary_hld} gives H\"older control of the horizontal gradient on non-characteristic boundary portions. In the proof of the centered Kor\'anyi-ball asymptotics this estimate will be applied away from the characteristic set, while the behavior near the characteristic points of the Kor\'anyi ball will be treated separately through the estimates recalled in Lemma~\ref{lemma:reg-ueps}.
}
\end{rem}

\subsection{The Pohozaev identity}

We denote by~$\mathdutchcal{D}$ the infinitesimal generator of the one-parameter group of non-isotropic dilations~\eqref{def_philambda}, namely
\begin{equation}\label{dilation_vector_field}
\mathdutchcal{D}
:=
\sum_{j=1}^{2n}x_j\partial_{x_j}
+
2t\partial_t .
\end{equation}
The following Pohozaev identity is one of the key tools in the paper.
\begin{lemma}[Pohozaev identity; see~\cite{GV00,GL92}]\label{pohozaev}
Let~$\Omega$ be a bounded~$C^1$-domain and let~$u\in C^2(\overline\Omega)$ solve
\[
-\Delta_Xu=f(u)
\qquad \text{in } \Omega,
\]
where~$f\in C(\r)$ and~$f(0)=0$. Set
\[
F(s):=\int_0^s f(\tau)\,{\rm d}\tau .
\]
Then,
\begin{eqnarray*}
&&\int_\Omega
\big(2QF(u)-(Q-2)u f(u)\big)\,{\rm d}\xi\\
&&\quad =
2\sum_{j=1}^{2n}
\int_{\partial\Omega}
\mathdutchcal{D} u\,X_j u\,\langle X_j,{n}_{\partial\Om}\rangle
\,{\rm d}\mathcal H^{Q-2}
-
\int_{\partial\Omega}
|\Dh u|^2\langle\mathdutchcal{D},{n}_{\partial\Om}\rangle
\,{\rm d}\mathcal H^{Q-2}
\\
&&\qquad
+
2\int_{\partial\Omega}
F(u)\langle\mathdutchcal{D},{n}_{\partial\Om}\rangle
\,{\rm d}\mathcal H^{Q-2}
+
(Q-2)\sum_{j=1}^{2n}
\int_{\partial\Omega}
u\,X_j u\,\langle X_j,{n}_{\partial\Om}\rangle
\,{\rm d}\mathcal H^{Q-2}\,,
\end{eqnarray*}
where~${n}_{\partial\Om}$ denotes the Euclidean exterior unit normal to~$\partial\Omega$ and ${\rm d}\mathcal H^{Q-2}$ denotes the Euclidean surface measure on~$\partial\Omega$.
\end{lemma}

\begin{rem}\label{rem:pohozaev-regularity}
{\rm
The identity above is first justified for smooth solutions and smooth domains. In the applications below, it will be applied through the standard approximation procedure allowed by the boundary regularity estimates available away from the characteristic set and by the additional control near the characteristic set.
}
\end{rem}

\subsection{Jerison-Lee extremals}

We conclude this preliminary section by recalling the explicit expression for the extremals in the sharp Folland-Stein-Sobolev inequality. The following result is due to Jerison and Lee~\cite{JL88}.

\begin{theorem}[Corollary~C in~\cite{JL88}]\label{thm_optimal}
Let~$2^\ast=2Q/(Q-2)$. Let~$U_{\rm JL}$ be given by
\begin{equation}\label{talentiane_2}
U_{\rm JL}(\xi)
=
c_0
\left(
(1+|z|^2)^2+t^2
\right)^{-\frac{Q-2}{4}}
\quad \xi=(z,t)\in\h,
\end{equation}
where the constant~$c_0>0$ is chosen according to the normalization
\[
-\Delta_XU_{\rm JL}=U_{\rm JL}^{2^\ast-1}
\quad
\text{in}~\h .
\]
Then, for every~$\lambda>0$ and every~$\xi_{\rm o}\in\h$, the function
\begin{equation}\label{eq:JL-family}
U_{{\rm JL},\lambda,\xi_{\rm o}}(\xi)
:=
\lambda^{-\frac{Q-2}{2}}
U_{\rm JL}\left(
\delta_{\lambda^{-1}}(\xi_{\rm o}^{-1}\circ\xi)
\right)
\end{equation}
is also a positive solution of
\[
-\Delta_XU_{{\rm JL},\lambda,\xi_{\rm o}}
=
U_{{\rm JL},\lambda,\xi_{\rm o}}^{2^\ast-1}
\quad
\text{in}~\h ,
\]
and realizes equality in the sharp Folland-Stein-Sobolev inequality~\eqref{folland}; namely,
$$
\|U_{{\rm JL},\lambda,\xi_{\rm o}}\|_{L^{2^\ast}(\h)}^{2^\ast}
=
\Sob
\|\Dh U_{{\rm JL},\lambda,\xi_{\rm o}}\|_{L^2(\h)}^{2^\ast}.
$$
Moreover, up to the invariances of the equation, all positive solutions of
\[
-\Delta_XV=V^{2^\ast-1}
\quad\text{in}~\h
\]
realizing equality are of this form. Consequently, all positive extremals for the homogeneous inequality~\eqref{folland} are precisely the positive multiples of the functions in the family~\eqref{eq:JL-family}.
\end{theorem}
\begin{rem}\label{rem:JL-normalization}
{\rm
The explicit profile in~\eqref{talentiane_2} is not, by itself, the normalization selected by the blow-up procedure. In the proof of Theorem~\ref{han}, the rescaled functions satisfy
\[
v_\eps(0)=\max v_\eps=1.
\]
Therefore the limiting Jerison-Lee extremal is the unique member of the family~\eqref{eq:JL-family}, centered at the identity, whose maximum is equal to one. More precisely, we fix once and for all
\begin{equation}\label{eq:def-blow-up-normalized-U}
U
:=
U_{{\rm JL},\lambda_\ast,0},
\qquad
\lambda_\ast
:=
U_{\rm JL}(0)^{\frac{2}{Q-2}} .
\end{equation}
Then

\[
-\Delta_XU=U^{2^\ast-1}
\quad\text{in}~\h 
\quad \text{and} \quad 
U(0)=\max_{\h}U=1 .
\]

Throughout the paper the constant appearing in the Green profile is defined by
\begin{equation}\label{eq:def-barc-n}
\bar c_n
:=
\int_{\h}U^{2^\ast-1}\,{\rm d}\xi ,
\end{equation}
with this blow-up normalization. With this convention, in the Green-profile limit proved below in the centered Kor\'anyi ball, one has
$$
\|u_\eps\|_{L^\infty(B_R(0))}u_\eps
\ \to \ 
\bar c_n G_{B_R(0)}(\cdot;0)
$$
locally uniformly in $B_R(0)\setminus\{0\}$, after the concentration point has been identified with the center. Thus~$\bar c_n$ is attached to the normalized extremal selected by~$v_\eps(0)=1$, not to an arbitrary algebraic representative of the Jerison-Lee family.

}
\end{rem}

\section{Asymptotic control via the Jerison-Lee extremals}\label{sec_asymptotic}

\subsection{Regularity properties for subelliptic equations}\label{sec_regprop}

We collect two local estimates for weak solutions of subelliptic equations in~$\h$ beginning with a local gain of integrability. In the two statements below the estimates are local in arbitrary balls whose centers do not need to belong to~$\overline\Omega$, since we adopt the convention that functions in~$\mathring{\mathcal{D}}^{1,2}_X(\Omega)$ are extended by zero outside~$\Omega$ whenever Sobolev inequalities are applied on balls intersecting the boundary.

\begin{lemma}\label{lemma:gain}
Let~$\Omega\subsetneq\h$ be a domain, let~$2<p_{\rm o}\leq p\leq 2^\ast$ and let~$u\in \mathring{\mathcal{D}}^{1,2}_X(\Omega)$ be a nonnegative weak solution of
\begin{equation}\label{pbm:han}
-\Delta_Xu=a(\xi)u^{p-1}
\quad \text{in}~\Omega,
\end{equation}
where~$a\in L^\infty(\Omega)$. Then, for every~$1<q\leq 2^\ast-1$, there exist constants~$\upsilon_{\rm o}>0$ and~$r_{\rm o}>0$, depending only on~$n$,~$p_{\rm o}$,~$q$, and~$\|a\|_{L^\infty(\Omega)}$, with the following property: if~$\eta\in\h$,~$0<r\leq r_{\rm o}$, and
\begin{equation}\label{eq:first-gain}
\left(\int_{\Omega\cap B_{2r}(\eta)}u^p\,{\rm d}\xi\right)^\frac{1}{p}
\leq
\upsilon_{\rm o},
\end{equation}
then
\begin{equation}\label{gain_int}
\left(\int_{\Om \cap B_r(\eta)}u^{\frac{(q+1)2^\ast}{2}}\,{\rm d}\xi\right)^\frac{2}{2^*(q+1)}
\leq
\frac{C}{r^{\frac{2}{q+1}}}\left(\int_{\Om \cap B_{2r}(\eta)}u^{q+1}\,{\rm d}\xi\right)^\frac{1}{q+1},
\end{equation}
where~$C>0$ depends only on~$n$ and~$q$.
\end{lemma}

\begin{proof}
 Let us assume~$\Om \cap B_{2r}(\eta) \neq \emptyset$, otherwise there is nothing to prove. Let~$\psi\in C^\infty_0(B_{2r}(\eta))$ be such that~$0\leq\psi\leq1$,~$\psi\equiv1$ in~$B_r(\eta)$ and~$|\Dh\psi|\leq C/r$, with~$C>0$ universal; see, e.~\!g.,~\cite[Lemma~3.6]{CGL93} for the existence of such a cut-off function. For any~$k\in \mathbb{N}$ let us define~$u_k:=\min\{u,k\}$. Since~$u\in \mathring{\mathcal{D}}^{1,2}_X(\Omega)$, the function~$\psi^2u u_k^{q-1}$ is an admissible test function in~\eqref{pbm:han} so that
\[
\int_{\Omega}
\langle \Dh u,\Dh(\psi^2u u_k^{q-1})\rangle\,{\rm d}\xi
=
\int_{\Omega}a(\xi)\psi^2u^{p}u_k^{q-1}\,{\rm d}\xi .
\]
On the other hand, by weak chain rule we have that
\begin{eqnarray*}
\int_{\Omega}
\langle \Dh u,\Dh(\psi^2u u_k^{q-1})\rangle\,{\rm d}\xi
& = & 2\int_{\Om \cap \{u < k\}} \psi u^{q}\langle \Dh u,\Dh\psi\rangle\,{\rm d}\xi + q\int_{\Om \cap \{u < k\}} \psi^2 u^{q-1}|\Dh u|^2\,{\rm d}\xi\notag\\
&& + k^{q-1}\int_{\Om \cap \{u \geq k\}} \left(2\psi u\langle \Dh u, \Dh \psi\rangle + \psi^2|\Dh u|^2\right)\,{\rm d}\xi\notag\\
 & \geq &
 \frac{q}{2}
 \int_{\Omega \cap \{u < k\}}\psi^2u^{q-1}|\Dh u|^2\,{\rm d}\xi+ \frac{k^{q-1}}{2}
 \int_{\Omega \cap \{u \geq k\}}\psi^2|\Dh u|^2\,{\rm d}\xi\\*[0.5ex]
 && - C_q
 \int_{\Omega}u^2u^{q-1}_k|\Dh\psi|^2\,{\rm d}\xi.
\end{eqnarray*}
Then, it follows that
\begin{multline}\label{eq:gain-caccioppoli}
\int_{\Omega \cap \{u < k\}}\psi^2u^{q-1}|\Dh u|^2\,{\rm d}\xi + k^{q-1}
 \int_{\Omega \cap \{u \geq k\}}\psi^2|\Dh u|^2\,{\rm d}\xi\\
\leq
C_q
\int_{\Omega}|a(\xi)|\psi^2u^{p}u_k^{q-1}\,{\rm d}\xi
+
C_q
\int_{\Omega}u^2u^{q-1}_k|\Dh\psi|^2\,{\rm d}\xi .
\end{multline}
Set now~$w_k:=\psi u^{\frac{q-1}{2}}_k u$ and note that since
\begin{eqnarray*}
 && |\Dh w_k|^2 \leq C_q u^{q+1}|\Dh \psi|^2 + C_q\psi^2u^{q-1}|\Dh u|^2 \quad \text{on}~\Om   \cap \{u < k\}\\*[1ex]
 \text{and} \quad && |\Dh w_k|^2 \leq C_q k^{q-1}u^2|\Dh \psi|^2 + C_qk^{q-1}\psi^2|\Dh u|^2 \quad \text{on}~\Om \cap \{u \geq  k\}\,,
\end{eqnarray*}
 using~\eqref{eq:gain-caccioppoli} and the Folland-Stein-Sobolev inequality~\eqref{folland}, we get
\begin{eqnarray}\label{eq:gain-sobolev}
(\Sob)^{-\frac{2}{2^\ast}}
\left( \int_{\Omega\cap B_{2r}(\eta)}
w_k^{2^\ast}\,{\rm d}\xi
\right)^{\frac{2}{2^\ast}}
& \leq & 
C_q \int_{\Omega\cap B_{2r}(\eta)}  |a(\xi)|\psi^2u^{p}u_k^{q-1}\,{\rm d}\xi\notag\\
&& + C_q\int_{\Omega\cap B_{2r}(\eta)}
u^2u^{q-1}_k|\Dh\psi|^2\,{\rm d}\xi .
\end{eqnarray}

For the first term on the right-hand side,  H\"older's Inequality gives
\begin{eqnarray*}
\int_{\Omega\cap B_{2r}(\eta)}
|a(\xi)|\psi^2u^{p}u_k^{q-1}\,{\rm d}\xi
& \leq & 
\|a\|_{L^\infty(\Omega)}
\int_{\Omega\cap B_{2r}(\eta)}
\big(\psi u^{\frac{q-1}{2}}_ku\big)^2 u^{p-2}\,{\rm d}\xi\\
&\leq &
\|a\|_{L^\infty(\Omega)}
\left(
\int_{\Omega\cap B_{2r}(\eta)}
 w_k^{2^\ast}\,{\rm d}\xi
\right)^{\frac{2}{2^\ast}}
\\
&& \times
\left(
\int_{\Omega\cap B_{2r}(\eta)}
u^{\frac{Q(p-2)}{2}}\,{\rm d}\xi
\right)^{\frac{2}{Q}} .
\end{eqnarray*}

Since~$p\leq 2^\ast$, one has~$Q(p-2)/2\leq p$.

Hence, by H\"older's Inequality
we get
$$
\left(
\int_{\Omega\cap B_{2r}(\eta)} u^{\frac{Q(p-2)}{2}}\,{\rm d}\xi \right)^{\frac{2}{Q}} 
\leq \left( \int_{\Omega\cap B_{2r}(\eta)}
u^p\,{\rm d}\xi \right)^{\frac{p-2}{p}} |\Omega\cap B_{2r}(\eta)|^{\frac{2}{Q}-\frac{p-2}{p}} .
$$
Choosing~$\upsilon_{\rm o}$ and~$r_{\rm o}$ sufficiently small, with constants depending only on~$n$,~$p_{\rm o}$,~$q$, and~$\|a\|_{L^\infty(\Omega)}$, we can absorb the corresponding term in the left-hand side of~\eqref{eq:gain-sobolev}. Therefore, by passing to the limit as~$k \to +\infty$, we get
$$
\left(\int_{\Omega\cap B_r(\eta)}u^{\frac{(q+1)2^\ast}{2}}\,{\rm d}\xi\right)^{\frac{2}{2^\ast}}
\leq
C r^{-2}\int_{\Omega\cap B_{2r}(\eta)}u^{q+1}\,{\rm d}\xi .
$$
Taking the~$(q+1)$-root on both sides gives~\eqref{gain_int}.
\end{proof}

We shall also use the following local boundedness estimate. It is a standard Moser iteration estimate for subelliptic equations with potential in~$L^\frac{q}{2}$,~$q>Q$.

\begin{lemma}\label{lemma:sup}
Let~$\Omega\subsetneq\h$ be a domain and let~$u\in \mathring{\mathcal{D}}^{1,2}_X(\Omega)$ be a nonnegative weak solution of
\begin{equation}\label{eq:moser-pbm}
-\Delta_Xu=h(\xi)u
\quad \text{in}~\Omega.
\end{equation}
Then, for every~$\eta\in\h$ and every~$r>0$ {such that~$h\in L^\frac{q}{2}(\Omega \cap B_{2r}(\eta))$, for some~$q>Q$,}
it holds
\begin{equation}\label{sup_est}
\sup_{\Omega\cap B_r(\eta)}u
\leq
C
\left( \frac{1}{|B_{2r}(\eta)|}
\int_{\Omega\cap B_{2r}(\eta)} u^{2^\ast}\,{\rm d}\xi \right)^{\frac{1}{2^\ast}},
\end{equation}
where~$C$ depends only on~$n$,~$q$, and on
$$
r^{2-\frac{2Q}{q}}
\|h\|_{L^\frac{q}{2}(\Omega\cap B_{2r}(\eta))}.
$$
In particular, if~$r$ is bounded from above and~$\|h\|_{L^\frac{q}{2}(\Omega)}$ is bounded, then~$C$ is uniform.
\end{lemma}
\begin{proof}
By the invariance of the equation, we let~$
\Omega_r:=\delta_\frac{1}{r}(\eta^{-1}\circ\Omega)$ and~$
h_r(\xi):=r^2h(\eta\circ\delta_r(\xi))$ so that~$
u_r(\xi):=u(\eta\circ\delta_r(\xi))$ solves~$
-\Delta_Xu_r=h_r(\xi)u_r$ in~$\Omega_r$. 
Thus, it is enough to prove the estimate for~$r=1$ and~$\eta=0$, with~$\Omega$ replaced by the rescaled domain~$\Omega_r$  and~$h$ replaced by~$h_r$. Moreover, let us assume that~$\Om \cap B_2(0) \neq \emptyset$. 

Let~$1\leq \sigma<\rho\leq2$, and choose~$\psi\in C^\infty_0(B_\rho(0))$ such that~$
0\leq\psi\leq1$,~$
\psi\equiv1$ in~$B_\sigma(0)$ and~$
|\Dh\psi|\leq C/(\rho-\sigma)$. For~$\theta>2$, we test~\eqref{eq:moser-pbm} with~$\psi^2uu^{\theta-2}_k$ where, as before,~$u_k:=\min\{u,k\}$. As in the proof of Lemma~\ref{lemma:gain}, via Young's Inequality and the Folland-Stein-Sobolev inequality for~$w_k := \psi u^{(\theta-2)/2}_ku$
we obtain
\begin{equation}\label{eq:moser-sobolev}
\|{w_k}\|_{L^{2^\ast}(\Omega)}^2
\leq
C\theta^2
\left(
\int_{\Omega}|h|\psi^2{u^{\theta-2}_ku^2}\,{\rm d}\xi
+
\int_{\Omega}{{u^2}u^{\theta-2}_k}|\Dh\psi|^2\,{\rm d}\xi
\right).
\end{equation}
By H\"older's Inequality,
$$
\int_{\Omega}|h|\psi^2{u^{\theta-2}_ku^2}\,{\rm d}\xi
\leq
\|h\|_{L^\frac{q}{2}(\Omega\cap B_\rho)}
\|{w_k}\|_{L^{\frac{2q}{q-2}}(\Omega)}^2 .
$$
Since~$q>Q$, we have that~$2q/(q-2) \in (2,2^*)$. 

Therefore, by interpolation between~$L^2$ and~$L^{2^\ast}$, for every~$\delta>0$,
\[
\|{w_k}\|_{L^{\frac{2q}{q-2}}(\Omega)}^2
\leq \delta\|{w_k}\|_{L^{2^\ast}(\Omega)}^2 + C\delta^{-\frac{Q}{q-Q}} \|{w_k}\|_{L^2(\Omega)}^2 .
\]

Choosing~{$\delta:= 1/(2C\theta^2(1+\|h\|_{L^\frac{q}{2}(\Omega\cap B_2)}))$}, 
we absorb the~$L^{2^*}$-norm of~{$w_k$} on the left-hand side
of~\eqref{eq:moser-sobolev}  
and after passing to the limit as~$k \to +\infty$, we obtain
\begin{equation}\label{eq:moser-recursive}
\|u\|_{L^{\theta\chi}(\Omega\cap B_\sigma)}
\leq
\left(
\frac{
C\theta^\alpha
\big(1+\|h\|_{L^\frac{q}{2}(\Omega\cap B_2)}^\alpha\big)
}
{(\rho-\sigma)^2}
\right)^{\frac{1}{\theta}}
\|u\|_{L^\theta(\Omega\cap B_\rho)},
\end{equation}
where~$
\chi:=\frac{2^\ast}{2}=\frac{Q}{Q-2}>1$,
and~$\alpha>0$ depends only on~$q$ and~$Q$.

We now proceed with a standard Moser's Iteration scheme in order to obtain the desired $L^\infty$-estimate. For this, for any~$j \in \mathbb{N}$ consider~$
\theta_j:=2^\ast\chi^j$ and~$\rho_j:=1+2^{-j}$, so that iterating~\eqref{eq:moser-recursive} with~$\theta=\theta_j$,~$\sigma=\rho_{j+1}$, and~$\rho=\rho_j$, via the convergence of the {series}
$$
\sum_{j=0}^\infty \frac{1}{\theta_j}<\infty \quad \text{and} \quad
\sum_{j=0}^\infty \frac{j}{\theta_j}<\infty,
$$
yields
$$
\sup_{\Omega\cap B_1}u
\leq C \|u\|_{L^{2^\ast}(\Omega\cap B_2)}.
$$
After scaling back to a general ball~$B_{2r}(\eta)$, this gives~\eqref{sup_est}.
\end{proof}

\subsection{\texorpdfstring{$\mathbb H$}{H}-Kelvin transform}\label{sec_kelvin}

We briefly recall the $\mathbb H$-Kelvin transform, which will be used in the proof of Theorem~\ref{han}. 

\begin{defn}\label{def_kelvin}
For every~$\xi=(z,t)\in\h\setminus\{0\}$, we define the $\mathbb H$-inversion by
\[
\begin{aligned}
  \kappa_{\mathbb H}& :\h\setminus\{0\}\longrightarrow\h\setminus\{0\},  \\ 
  \kappa_{\mathbb H}(z,t)
& :=
\left(
-\frac{z}{|z|^2-it},
-\frac{t}{|z|^4+t^2}
\right).
\end{aligned}
\]
Given a function~$u:\h\to\r$, its $\mathbb H$-Kelvin transform is the function~$u^\sharp$ defined by
\[
u^\sharp(\xi)
:=
|\xi|_{\mathbb H}^{-(Q-2)}
u\big(\kappa_{\mathbb H}(\xi)\big)
\quad
\xi\in\h\setminus\{0\}.
\]
\end{defn}

With this convention one has
\begin{equation}\label{H-inve-prop}
\kappa_{\mathbb H}(\kappa_{\mathbb H}(\xi))=\xi 
\quad {\text{and}} \quad
|\kappa_{\mathbb H}(\xi)|_{\mathbb H}
=
|\xi|_{\mathbb H}^{-1}
\quad
\forall \xi\neq0.
\end{equation}
{
\begin{rem}{\rm
    We shall denote with~$\Om^\sharp$ the image of a generic domain~$\Om$ under the $\mathbb{H}$-inversion. In the case when~$\Om$ is a neighborhood of~$\infty$, i.~\!e. when there exists~$B_R(0)$ such that~$(\h \setminus B_R(0)) \subsetneq \Om$, then~$\Om^\sharp$ is a punctured neighborhood of the origin~$0$, i.~\!e.~$\Om^\sharp = A \setminus \{0\}$ for some open set~$A$ containing the group identity.}
\end{rem}
}

We shall use the following standard properties of the Kelvin transform.

\begin{prop}[Theorem~2.3.5 and Lemma~2.3.6 in~\cite{IV11}]\label{prop:H-Kelv-iso-CR-lap}
The following properties hold:
\begin{enumerate}[\it (i)]
    \item  {Let~$\Om^\sharp$ be the image of~$\Om$ with respect to the $\mathbb{H}$-inversion. Then, the Kelvin transform}
is an isometry between~$\mathring{\mathcal{D}}^{1,2}_X(\Omega)$ and~$\mathring{\mathcal{D}}^{1,2}_X(\Omega^\sharp)$. Moreover, for every measurable set~$E\subset\Omega$,
\begin{equation}\label{eq:kelvin-critical-measure}
\int_{\kappa_{\mathbb H}(E\setminus\{0\})}
|u^\sharp(\xi)|^{2^\ast}\,{\rm d}\xi
=
\int_E |u(\xi)|^{2^\ast}\,{\rm d}\xi,
\end{equation}
whenever one of the two integrals is finite.
     \item For~$p>0$, let~$u$ be a solution of
\[
\begin{cases}
-\Delta_Xu=u^p
\quad
\text{in}~\Omega\\
u \in \mathring{\mathcal{D}}^{1,2}_X(\Omega), \quad u \geq 0.
\end{cases}
\]
Then, its $\mathbb H$-Kelvin transform~$u^\sharp$ satisfies
\[
\begin{cases}
-\Delta_Xu^\sharp(\xi)
=
|\xi|_{\mathbb H}^{p(Q-2)-(Q+2)}
\big(u^\sharp(\xi)\big)^p
\quad
\text{in}~\Omega^\sharp,\\
u^\sharp\in \mathring{\mathcal{D}}^{1,2}_X(\Omega^\sharp),
\quad
u^\sharp\geq0.
\end{cases}
\]
 In particular, when~$p=2^\ast-1$, the critical equation is invariant under the $\mathbb H$-Kelvin transform
\end{enumerate}
\end{prop}

\begin{rem}\label{rem:kelvin-subcritical}
{\rm
For the subcritical exponent $p=2^\ast-\eps-1$, the transformed equation becomes
$$
-\Delta_Xu^\sharp(\xi)
=
|\xi|_{\mathbb H}^{-\eps(Q-2)}
\big(u^\sharp(\xi)\big)^{2^\ast-\eps-1}.
$$
This weight is locally bounded away from the origin and is the only change with respect to the critical case. This observation will be used below when applying local estimates to Kelvin-transformed functions.
}
\end{rem}

\subsection{Convergence of the blow-up}\label{subsec:boundary-exclusion}

We prove here the convergence of the blow-up of the subcritical maximizer~$u_\eps$. This is done by combining asymptotic sharpness of the critical quotient with the Brezis--Lieb lemma.

\begin{lemma}\label{lemma:tail-compactness-v}
Let~$u_\eps$ be the scalar-normalized maximizers solving
\[
\begin{cases}
    -\Delta_Xu_\eps=u_\eps^{2^\ast-\eps-1} & \text{in}~\Om,\\
    u_\eps =0 & \text{on}~\partial\Om.
\end{cases}
\]
Let~$M_\eps:=\|u_\eps\|_{L^\infty(\Omega)}$, and let~$\eta_\eps\in\Omega$ be such that $u_\eps(\eta_\eps)=M_\eps$. Set $\rho_\eps:=M_\eps^{-({2^\ast-\eps-2})/{2}}$,
and define
$$
v_\eps(\zeta)
:= M_\eps^{-1}u_\eps(\eta_\eps\circ\delta_{\rho_\eps}(\zeta)),
\qquad \zeta\in \Omega_{\eta_\eps,\rho_\eps}:=\delta_{\frac{1}{\rho_\eps}}(\eta_\eps^{-1}\circ\Omega)
$$
Assume that, up to subsequences, the following assumptions do hold:
\begin{enumerate}[\it (i)]
    \item $1\leq M_\eps^\eps\leq C $, for some~$C>1$ independent of~$\eps$;
    \item $ v_\eps \to U $ locally uniformly in~$\h$ and locally in~$\mathring{\mathcal{D}}^{1,2}_X(\h)$, where~$U$ is the blow-up normalized Jerison-Lee extremal fixed in~\eqref{eq:def-blow-up-normalized-U}.
\end{enumerate}
Then, after extending~$v_\eps$ by zero outside~$\Omega_{\eta_\eps,\rho_\eps}$, one has
\begin{equation}\label{eq:strong-L2star-v}
v_\eps\to U
\quad
\text{in}~L^{2^\ast}(\h)\,,
\end{equation}
and
\begin{equation}\label{eq:uniform-tail-v}
\lim_{R\to+\infty}
\limsup_{\eps\to0}
\int_{\Omega_{\eta_\eps,\rho_\eps}\setminus B_R(0)}
v_\eps^{2^\ast}\,{\rm d}\zeta
=0.
\end{equation}
\end{lemma}

\begin{proof}
We extend~$v_\eps$ by zero outside~$\Omega_{\eta_\eps,\rho_\eps}$. By the change of variables defining~$v_\eps$, one has
\begin{equation} \label{eq:veps-grad-scaling}
\|\Dh v_\eps\|_{L^2(\h)}^2 =M_\eps^{-\frac{\eps(Q-2)}{2}} \|\Dh u_\eps\|_{L^2(\Omega)}^2 .
\end{equation}
By the scalar normalization of the maximizers, equivalently by the definition of~$(\Sob)^{-\frac{Q-2}{2}}$ in~\eqref{eq:def-Sstar-normalized}, the right-hand side is bounded uniformly in~$\eps$. Hence, using also the bound in~{\it (i)} gives
$$
\sup_{\eps>0} \|\Dh v_\eps\|_{L^2(\h)} <+\infty .
$$
Thus the zero extensions of~$v_\eps$ are bounded in~$\mathring{\mathcal D}^{1,2}_X(\h)$.

We next record the asymptotic sharpness of the critical quotient. For every nonzero~$w\in\mathring{\mathcal D}^{1,2}_X(\h)$, set
$$
\mathcal Q(w)
:=
\frac{\displaystyle\int_{\h}|w|^{2^\ast}\,{\rm d}\xi}
{\displaystyle\left(\int_{\h}|\Dh w|^2\,{\rm d}\xi\right)^{2^\ast/2}}.
$$
Since the critical quotient is invariant under multiplication by positive constants, left translations, and anisotropic dilations,
$$
\mathcal Q(v_\eps)=\mathcal Q(u_\eps).
$$
Let
$$
E_\eps
:=
\int_\Omega|\Dh u_\eps|^2\,{\rm d}\xi
=
\int_\Omega u_\eps^{2^\ast-\eps}\,{\rm d}\xi.
$$
By H\"older's inequality,
$$
\int_\Omega u_\eps^{2^\ast}\,{\rm d}\xi
\geq
|\Omega|^{1-\frac{2^\ast}{2^\ast-\eps}}
E_\eps^{\frac{2^\ast}{2^\ast-\eps}}.
$$
Therefore
$$
\mathcal Q(u_\eps)
\geq
|\Omega|^{1-\frac{2^\ast}{2^\ast-\eps}}
E_\eps^{\frac{2^\ast}{2^\ast-\eps}-\frac{2^\ast}{2}}
\to
\Sob,
$$
where~\eqref{eq:def-Sstar-normalized} has been used. Since~$\mathcal Q(u_\eps)\leq\Sob$, we obtain
\begin{equation}\label{eq:sharp-quotient-v}
\lim_{\eps\to0^+}\mathcal Q(v_\eps)=\Sob.
\end{equation}

Set
$$
r_\eps:=v_\eps-U.
$$
Since~$\{v_\eps\}$ is bounded in~$\mathring{\mathcal D}^{1,2}_X(\h)$ and~$v_\eps\to U$ locally uniformly, after passing to a subsequence we have
$$
v_\eps\rightharpoonup U
\quad\text{in }\mathring{\mathcal D}^{1,2}_X(\h),
\qquad
v_\eps\to U
\quad\text{a.e. in }\h.
$$
Hence
\begin{equation}\label{eq:tail-energy-splitting}
\|\Dh v_\eps\|_{L^2(\h)}^2
=
\|\Dh U\|_{L^2(\h)}^2
+
\|\Dh r_\eps\|_{L^2(\h)}^2
+o_\eps(1).
\end{equation}
By the Brezis--Lieb lemma,
\begin{equation}\label{eq:tail-critical-splitting}
\|v_\eps\|_{L^{2^\ast}(\h)}^{2^\ast}
=
\|U\|_{L^{2^\ast}(\h)}^{2^\ast}
+
\|r_\eps\|_{L^{2^\ast}(\h)}^{2^\ast}
+o_\eps(1).
\end{equation}
Set
$$
A_0:=\|\Dh U\|_{L^2(\h)}^2>0,
\qquad
B_\eps:=\|\Dh r_\eps\|_{L^2(\h)}^2.
$$
Since~$U$ realizes equality in~\eqref{folland},
$$
\|U\|_{L^{2^\ast}(\h)}^{2^\ast}
=
\Sob A_0^{2^\ast/2},
$$
whereas~\eqref{folland} gives
$$
\|r_\eps\|_{L^{2^\ast}(\h)}^{2^\ast}
\leq
\Sob B_\eps^{2^\ast/2}.
$$
If, along a further subsequence,~$B_\eps\to B\geq0$, then~\eqref{eq:sharp-quotient-v},~\eqref{eq:tail-energy-splitting}, and~\eqref{eq:tail-critical-splitting} yield
$$
1
\leq
\frac{A_0^{2^\ast/2}+B^{2^\ast/2}}
{(A_0+B)^{2^\ast/2}}
\leq1.
$$
Since~$2^\ast/2>1$ and~$A_0>0$, the quotient in the middle is strictly smaller than~$1$ whenever~$B>0$. Thus~$B=0$. Since this holds for every convergent subsequence of~$\{B_\eps\}$,
$$
\|\Dh(v_\eps-U)\|_{L^2(\h)}\to0.
$$
The Sobolev inequality then gives~\eqref{eq:strong-L2star-v}.

 It remains to prove the tail estimate. Since~$U\in L^{2^\ast}(\h)$, for every~$\delta>0$ there exists~$R_\delta>0$ such that
$$
\int_{\h\setminus B_{R_\delta}(0)} U^{2^\ast}\,{\rm d}\zeta <\delta .
$$
Using~\eqref{eq:strong-L2star-v}
we get
$$
\begin{aligned}
\limsup_{\eps\to0} 
\int_{\Omega_{\eta_\eps,\rho_\eps}\setminus B_{R_\delta}(0)}
v_\eps^{2^\ast}\,{\rm d}\zeta
&\leq
C \limsup_{\eps\to0} \int_{\h}|v_\eps-U|^{2^\ast}\,{\rm d}\zeta + C \int_{\h\setminus B_{R_\delta}(0)} U^{2^\ast}\,{\rm d}\zeta  \\
&\leq C\delta .
\end{aligned}
$$
Letting~\(\delta\to0\) gives~\eqref{eq:uniform-tail-v}.
\end{proof}

\subsection{Proof of Theorem~\ref{han}}\label{proof_han}
Let~$u_\eps$ be the scalar-normalized maximizer introduced in~\eqref{eq:normalized-critical-subcritical}; namely,
\begin{equation}\label{eq:han-normalized-equation}
\begin{cases}
-\Delta_Xu_\eps=u_\eps^{2^\ast-\eps-1}
& \text{in}~\Om,\\
u_\eps=0
& \text{on}~\partial\Omega
\end{cases}
\end{equation}
and
\begin{equation}\label{eq_han1}
\int_\Omega u_\eps^{2^\ast -\eps}\,{\rm d}\xi
=
(\Sob)^{-\frac{Q-2}{2}}+{{\it o}_\eps}(1)
\quad
\text{as}~\eps\searrow0,
\end{equation}
where~${\it o}_\eps(1)$ denotes a quantity that goes to~$0$ when~$\eps\searrow 0$.

Moreover, let us assume there exist~$\delta>0$,~$C_0>0$, and~$\eps_{\rm o}>0$ such that
\begin{equation}\label{eq:away-boundary}
\sup_{\{\xi\in\Om:\ \operatorname{dist}_{\mathbb R^{2n+1}}(\xi,\partial\Om)<\delta\}}u_\eps(\xi)
\leq C_0
\qquad
\forall\eps\in(0,\eps_{\rm o}).
\end{equation}

We divide the proof into several steps.

\vspace{2mm}
{\it Step 1 The maximum of~$u_\eps$ diverges; i.~\!e., 
\begin{equation}\label{han_supremum_limit}
M_\eps:=\|u_\eps\|_{L^\infty(\Omega)} \to +\infty \quad \text{as}~\eps\searrow0.
\end{equation}
}
 Assume by contradiction that, along a sequence~$\eps_k\searrow0$, the functions
$u_{\eps_k}$ are uniformly bounded in~$L^\infty(\Omega)$. Testing
\eqref{eq:han-normalized-equation} with~$u_{\eps_k}$ and using~\eqref{eq_han1}, we get
$$
\int_\Omega |\Dh u_{\eps_k}|^2\,{\rm d}\xi = \int_\Omega u_{\eps_k}^{2^\ast-\eps_k}\,{\rm d}\xi = (\Sob)^{-\frac{Q-2}{2}}+o_{\eps_k}(1).
$$
Hence $\{u_{\eps_k}\}_{k\in\mathbb N}$ is bounded in
$\mathring{\mathcal D}^{1,2}_X(\Omega)$. Up to a subsequence,
$$
u_{\eps_k}\rightharpoonup u \quad\text{weakly in } \mathring{\mathcal D}^{1,2}_X(\Omega),
$$
and, by~\cite[Theorem 3.2]{GL92}, 
$$
u_{\eps_k}\to u \quad\text{a.~\!e. in }\Omega \quad\text{and strongly in }L^q(\Omega) \quad\text{for every }1\le q<2^\ast .
$$

Therefore, dominated convergence and~\eqref{eq_han1} give
$$
\int_\Omega u^{2^\ast}\,{\rm d}\xi
=
(\Sob)^{-\frac{Q-2}{2}}>0,
$$

so that~$u\not\equiv0$. Moreover, by weak lower semicontinuity,
$$
\int_\Omega |\Dh u|^2\,{\rm d}\xi
\leq
\liminf_{k\to+\infty}
\int_\Omega |\Dh u_{\eps_k}|^2\,{\rm d}\xi
=
(\Sob)^{-\frac{Q-2}{2}} .
$$
Hence
$$
\mathcal Q(u)
=
\frac{\displaystyle\int_\Omega u^{2^\ast}\,{\rm d}\xi}
{\displaystyle\left(\int_\Omega |\Dh u|^2\,{\rm d}\xi\right)^{2^\ast/2}}
\ge
\frac{(\Sob)^{-\frac{Q-2}{2}}}{(\Sob)^{-\frac{Q}{2}}}
=
\Sob .
$$
By the sharp Sobolev inequality, $\mathcal Q(u)\leq\Sob$. Therefore $\mathcal Q(u)=\Sob$, and $u$ is a nontrivial maximizer for the critical problem~\eqref{critica} on the bounded domain~$\Omega$, contradicting the non-attainment of the critical Sobolev constant. Hence~\eqref{han_supremum_limit} follows.

{{Let~$\eta_\eps\in\Omega$ be such that
\begin{equation}\label{eq:def-M-eps}
u_\eps(\eta_\eps)=M_\eps 
.
\end{equation}
Set
\begin{eqnarray}
    \label{eq:def-rho-eps-proof}
\rho_\eps
& := &
M_\eps^{-\frac{2^\ast-\eps-2}{2}},\\
\Omega_{\eta_\eps,\rho_\eps}
& := &
\delta_{\rho_\eps^{-1}}\big(\eta_\eps^{-1}\circ\Omega\big), \label{eq:def-Omega-eps-han}\\ v_\eps(\zeta)
& := &
M_\eps^{-1} u_\eps\big(\eta_\eps\circ\delta_{\rho_\eps}(\zeta)\big), \quad \zeta\in{\Om_{\eta_\eps,\rho_\eps}}. \label{han_v_eps}
\end{eqnarray}}

{Moreover, the change of variables
$$
\xi=\eta_\eps\circ\delta_{\rho_\eps}(\zeta)
$$
together with~\eqref{eq:veps-grad-scaling} gives
\begin{equation}\label{eq:scaling-v-eps}
\int_{\Om_{\eta_\eps,\rho_\eps}}|\Dh v_\eps|^2\,{\rm d}\zeta = 
M_\eps^{-\frac{Q-2}{2}\eps}
\int_\Omega u_\eps^{2^\ast-\eps}\,{\rm d}\xi .
\end{equation}}}

Then
\begin{equation}\label{eq_han3}
\begin{cases}
-\Delta_Xv_\eps=v_\eps^{2^\ast-\eps-1}
\,\, \text{in}~{\Om_{\eta_\eps,\rho_\eps}},\\
v_\eps(0)=1,~0\leq v_\eps\leq1.
\end{cases}
\end{equation}

Since~$M_\eps\to+\infty$, in light of~\eqref{eq:away-boundary} it follows that, for~$\eps\in(0,\eps_{\rm o})$,
\begin{equation}\label{eq:eta-away-boundary-han}
\operatorname{dist}_{\mathbb R^{2n+1}}(\eta_\eps,\partial\Omega)\geq\delta.
\end{equation}
Since~$\rho_\eps\to0$, up to a subsequence,
\begin{equation}\label{eq:dist-over-rho}
\eta_\eps\to\bar\xi
\in
\{\xi\in\Omega:\operatorname{dist}_{\mathbb R^{2n+1}}(\xi,\partial\Omega)\geq\delta\},
\qquad
\frac{\operatorname{dist}_{\mathbb R^{2n+1}}(\eta_\eps,\partial\Omega)}{\rho_\eps}
\to+\infty.
\end{equation}
Thus~$\Omega_{\eta_\eps,\rho_\eps}$ exhausts the whole Heisenberg group in the sense of compact subsets.

Since $0\leq v_\eps\leq1$, local regularity estimates for subelliptic equations imply that $\{v_\eps\}$ is locally equicontinuous on every compact subset of $\h$. {Moreover, since~$M_\eps \geq 1$ by~\eqref{han_supremum_limit},~\eqref{eq:scaling-v-eps} and that~$u_\eps$ is a scalar-normalized maximizer, we get that~$v_\eps$ has uniformly bounded Sobolev energy. Therefore,} after extending~$v_\eps$ by zero outside~$\Omega_{\eta_\eps,\rho_\eps}$, up to a subsequence,
\begin{equation}\label{eq:v-eps-local-limit}
\begin{aligned}
  &  v_\eps \to  v_\infty \quad \text{locally uniformly in}~\h\,,\\
  \text{and} \quad & v_\eps \rightharpoonup v_\infty \quad \text{in}~\mathring{\mathcal{D}}^{1,2}_X(\h)\,,
\end{aligned}
\end{equation}
where
\begin{equation}\label{eq:han-limit-bubble1}
v_\infty(0)=1, \quad 0\leq v_\infty\leq1 .
\end{equation}

\vspace{2mm}
{\it Step 2  The quantity~$M_\eps^\eps$ is uniformly bounded, i.~\!e.
\begin{equation}\label{bound_M_eps}
 1 \leq M_\eps^\eps\leq C .
\end{equation}
}

{The lower bound is trivial and follows from~\eqref{han_supremum_limit}. Let us prove the upper bound.} Passing to the limit in~\eqref{eq_han3}, {via~\eqref{eq:v-eps-local-limit}}, we first obtain {that~$v_\infty$ is a weak solution to}
$$
-\Delta_Xv_\infty=v_\infty^{2^\ast-1} \quad \text{in}~\h.
$$

In particular~$v_\infty\not\equiv0$. 

{Fix}~$\sigma>0$ such that
\begin{equation}\label{importante}
\int_{B_\sigma(0)}v_\infty^{2^\ast}\,{\rm d}\zeta>0 .
\end{equation}
Since~{$\Om_{\eta_\eps,\rho_\eps}$} exhausts~$\h$,~$B_\sigma(0)$ is contained in~$\Om_{\eta_\eps,\rho_\eps}$ for~$\eps$ small enough, {and by~\eqref{eq:v-eps-local-limit} we get}
$$\lim_{\eps\to0}
\int_{B_\sigma(0)}v_\eps^{2^\ast-\eps}\,{\rm d}\zeta
= \int_{B_\sigma(0)}v_\infty^{2^\ast}\,{\rm d}\zeta >0 .
$$
On the other hand, {by~\eqref{eq:veps-grad-scaling} and~\eqref{eq:def-rho-eps-proof}}
we have
{ 
$$
\int_{B_\sigma(0)}v_\eps^{2^\ast-\eps}\,{\rm d}\zeta
\leq M_\eps^{-\frac{Q-2}{2}\eps}
\int_\Omega u_\eps^{2^\ast-\eps}\,{\rm d}\xi .
$$
}
{which, combined with}~\eqref{eq_han1} {and~\eqref{importante}}   gives the upper bound in~\eqref{bound_M_eps} {~$\eps>0$ sufficiently small.}

\vspace{2mm}
{\it Step 3 The maximum points converge to the concentration point; i.~\!e.
\begin{equation}\label{conv_xi_eps-proof}
\eta_\eps\to\xi_{\rm o}.
\end{equation}
}

By Theorem~\ref{cor_concentration}, applied before the scalar normalization and then rescaled as in~\eqref{eq:scalar-normalized-concentration-measure}, up to a subsequence there exists~$\xi_{\rm o}\in\overline\Omega$ such that
\begin{equation}\label{eq:han-proof-weakmea}
|\Dh u_\eps|^2\,{\rm d}\xi
\tows
(\Sob)^{-\frac{Q-2}{2}}\,\boldsymbol{\delta}_{\xi_{\rm o}}
\quad
\text{in}~\mea .
\end{equation}
We {are left to} show that~$\bar\xi=\xi_{\rm o}$ {with~$\bar{\xi}$ given in~\eqref{eq:dist-over-rho}. For this, let us fix}~$R>0$ be such that
$$
\int_{B_R(0)}|\Dh v_\infty|^2\,{\rm d}\zeta>0 .
$$
Such an~$R$ exists because~$v_\infty(0)=1$ and~$v_\infty\in\mathring{\mathcal D}^{1,2}_X(\h)$. Indeed, if
$$
\int_{B_R(0)}|\Dh v_\infty|^2\,{\rm d}\zeta=0
$$
for every~$R>0$, then~$\Dh v_\infty=0$ almost everywhere in~$\h$, so~$v_\infty$ is constant. This contradicts~$v_\infty(0)=1$ and~$v_\infty\in L^{2^\ast}(\h)$.

By lower semicontinuity, \eqref{eq:veps-grad-scaling} and~\eqref{eq:v-eps-local-limit} we get
$$
\begin{aligned}
       & 0< \int_{B_R(0)}|\Dh v_\infty|^2\,{\rm d}\zeta
\leq \liminf_{\eps\to0}
\int_{B_R(0)}|\Dh v_\eps|^2\,{\rm d}\zeta \\
\text{and} \quad & \int_{B_R(0)}|\Dh v_\eps|^2\,{\rm d}\zeta
= M_\eps^{-\frac{\eps(Q-2)}{2}}\int_{B_{R\rho_\eps}(\eta_\eps)}
|\Dh u_\eps|^2\,{\rm d}\xi.
\end{aligned}
$$
which gives
\[
\liminf_{\eps\to0}
\int_{B_{R\rho_\eps}(\eta_\eps)}
|\Dh u_\eps|^2\,{\rm d}\xi
>0 .
\]
Up to a standard regularization via cutoff, the preceding local energy lower bound, together with~\eqref{han_supremum_limit},~\eqref{eq:def-rho-eps-proof}, and~\eqref{eq:dist-over-rho}, shows that the limiting concentration measure assigns positive mass to every neighborhood of~$\bar\xi$. Hence~\eqref{eq:han-proof-weakmea} forces~$\bar\xi=\xi_{\rm o}$, which gives~\eqref{conv_xi_eps-proof}.

\vspace{2mm}

{\it Step 4 The rescaled functions converge to the normalized Jerison-Lee extremal~{$U$ given in~\eqref{eq:def-blow-up-normalized-U}}.}

\vs 
{Since~$v_\infty$ has finite horizontal energy by lower semicontinuity and~\eqref{eq:v-eps-local-limit} we get
\[
\int_{\h}|\Dh v_\infty|^2\,{\rm d}\zeta<  +\infty .
\]
}

{
Since $v_\infty\not\equiv0$ and $v_\infty\ge0$, the strong maximum principle gives
$$
v_\infty>0 \quad\text{in }\h .
$$

Moreover~$v_\infty$ has finite horizontal energy. Testing the limiting equation with~$v_\infty$, up to a standard regularization via cutoff, gives
$$
\int_{\h}|\Dh v_\infty|^2\,{\rm d}\zeta
=
\int_{\h}v_\infty^{2^\ast}\,{\rm d}\zeta .
$$
Set
$$
A:=\int_{\h}|\Dh v_\infty|^2\,{\rm d}\zeta .
$$
Since $v_\infty\not\equiv0$, the Sobolev inequality gives
$$
A\leq \Sob A^{2^\ast/2},
$$
and hence $A\geq(\Sob)^{-\frac{Q-2}{2}}$. On the other hand, by lower semicontinuity and~\eqref{eq:scaling-v-eps},
$$
A
\leq
\liminf_{\eps\to0}
\int_{\Om_{\eta_\eps,\rho_\eps}}|\Dh v_\eps|^2\,{\rm d}\zeta
\leq
\limsup_{\eps\to0}
\int_\Omega u_\eps^{2^\ast-\eps}\,{\rm d}\xi
=
(\Sob)^{-\frac{Q-2}{2}} ,
$$
because $M_\eps^{-\frac{(Q-2)\eps}{2}}\leq1$. Therefore $A=(\Sob)^{-\frac{Q-2}{2}}$ and equality holds in the sharp Folland--Stein inequality. By Theorem~\ref{thm_optimal}, $v_\infty$ belongs to the Jerison--Lee family. Finally, by~\eqref{eq:han-limit-bubble1},
$$
v_\infty(0)=\max_{\h}v_\infty=1 .
$$
By~\eqref{eq:def-blow-up-normalized-U}, this gives
$$
v_\infty=U .
$$
}

{{We show that the convergence is actually locally strong in~$\mathring{\mathcal D}^{1,2}_X$. Let $B_r(0)\Subset\h$ and choose
$\psi\in C^\infty_0(\h)$ such that $\psi\equiv1$ on~$B_r(0)$; {see, e.~\!g.,~~\cite[Lemma~3.6]{CGL93} for the existence of such a cut-off function}. Since
$\Omega_{\eta_\eps,\rho_\eps}$ exhausts~$\h$, we have
$\operatorname{supp}(\psi)\Subset\Omega_{\eta_\eps,\rho_\eps}$ for $\eps$ small enough. Set
$$
w_\eps:=v_\eps-U .
$$
Subtracting the equations solved by~$v_\eps$ and~$U$ and testing with
$\psi^2w_\eps$, we obtain
\begin{equation}\label{eq:han-bubble-strong-conv1}
\int_{\h}\psi^2|\Dh w_\eps|^2\,{\rm d}\zeta
=
\int_{\h}\psi^2
\big(v_\eps^{2^\ast-\eps-1}-U^{2^\ast-1}\big)w_\eps\,{\rm d}\zeta - 2\int_{\h}\psi w_\eps
\langle\Dh w_\eps,\Dh\psi\rangle\,{\rm d}\zeta .
\end{equation}
By Cauchy's Inequality,
$$
\left| 2\int_{\h}\psi w_\eps
\langle\Dh w_\eps,\Dh\psi\rangle\,{\rm d}\zeta \right|
\leq \frac{1}{2}
\int_{\h}\psi^2|\Dh w_\eps|^2\,{\rm d}\zeta + 2\|\Dh\psi\|_{L^\infty(\h)}^2
\int_{\operatorname{supp}(\psi)}w_\eps^2\,{\rm d}\zeta .
$$
Therefore, from~\eqref{eq:han-bubble-strong-conv1},
$$
\int_{B_r(0)}|\Dh w_\eps|^2\,{\rm d}\zeta
\leq 2\int_{\h}\psi^2
\left|v_\eps^{2^\ast-\eps-1}-U^{2^\ast-1}\right|
|w_\eps|\,{\rm d}\zeta
+ 4\|\Dh\psi\|_{L^\infty(\h)}^2
\int_{\operatorname{supp}(\psi)}w_\eps^2\,{\rm d}\zeta .
$$
Since $v_\eps\to U$ uniformly on $\operatorname{supp}(\psi)$ and
$0\le v_\eps\le1$, the right-hand side tends to zero. Hence
\[
\int_{B_r(0)}|\Dh(v_\eps-U)|^2\,{\rm d}\zeta\to0 .
\]}}
Thus, {Lemma~\ref{lemma:tail-compactness-v} we can upgrade~\eqref{eq:v-eps-local-limit} to
\begin{align}
  &  v_\eps \to U \quad \text{locally uniformly in}~\h\,,\label{eq:v-limit-U-unif}\\
  & v_\eps \to U \quad \text{locally in}~\mathring{\mathcal{D}}^{1,2}_X(\h)\,, \label{eq:v-limit-U-sob}\\
&   v_\eps\to U
\quad
\text{in}~L^{2^\ast}(\h)\,, \label{eq:strong-L2star-v-proof} \\
\text{and} \quad & \lim_{R\to+\infty}
\limsup_{\eps\to0}
\int_{\Omega_{\eta_\eps,\rho_\eps}\setminus B_R(0)}
v_\eps^{2^\ast}\,{\rm d}\zeta
=0 . \label{eq:uniform-tail-v-proof}
\end{align}
}

\vspace{2mm}
{\it Step 5 Kelvin transform estimate and global bubble bound.}
\vspace{1mm}

We now prove the global pointwise estimate
\begin{equation}\label{boun_max_seq2}
v_\eps(\zeta)\leq C U(\zeta)
\quad
\text{for every}~\zeta\in {\Omega_{\eta_\eps,\rho_\eps}},
\end{equation}
{with~$C=C(n,\Omega,\delta,C_0)>0$ being independent of~$\eps$.} Once~\eqref{boun_max_seq2} is established, the estimate~\eqref{bound_max_seq} follows immediately from the definition~\eqref{han_v_eps} of~$v_\eps$.

Let~$v_\eps^\sharp$ be the Kelvin transform of~$v_\eps$, and let
$$
\Omega_{\eta_\eps,\rho_\eps}^\sharp:=\kappa_{\mathbb H}(\Omega_{\eta_\eps,\rho_\eps}\setminus\{0\}).
$$
Since $\Omega$ is bounded and $\eta_\eps\to\xi_{\rm o}\in\Omega$, there exists a constant $C_\Omega>0$, independent of $\eps$, such that
\begin{equation}\label{eq:Omega-eps-exterior-radius}
\Omega_{\eta_\eps,\rho_\eps}
\subset
B_{R_\eps}(0)
\quad \text{and} \quad
R_\eps:=\frac{C_\Omega}{\rho_\eps}\,,
\end{equation}
{which yields that}
\begin{equation}\label{eq:Omega-sharp-away-origin}
\Omega_{\eta_\eps,\rho_\eps}^\sharp
\subset
\h\setminus B_{R_\eps^{-1}}(0).
\end{equation}
{By Proposition~\ref{prop:H-Kelv-iso-CR-lap}~{\it (ii)} and Remark~\ref{rem:kelvin-subcritical} give} that~$v_\eps^\sharp$ satisfies
\begin{equation}\label{eq:kelvin-v-eps}
-\Delta_Xv_\eps^\sharp(\zeta)
= |\zeta|_{\mathbb H}^{-\eps(Q-2)}
\big(v_\eps^\sharp(\zeta)\big)^{2^\ast-\eps-1}
\quad
\text{in}~\Omega_{\eta_\eps,\rho_\eps}^\sharp .
\end{equation}
Since~$0\leq v_\eps\leq1$, we obtain the following
\begin{equation}\label{eq:kelvin-rough-bound}
0\leq v_\eps^\sharp(\zeta)
\leq
|\zeta|_{\mathbb H}^{2-Q}\,;
\end{equation}
i.~\!e.,~$v_\eps^\sharp$ is uniformly bounded away from the origin. 

Hence it is enough first to prove a uniform bound for $v_\eps^\sharp$ in a fixed neighborhood of the origin. We shall show that there exist
$r_\sharp>0$ and $C>0$, independent of~$\eps$, such that
\begin{equation}\label{boun_max_seq3}
\sup_{\Omega_{\eta_\eps,\rho_\eps}^\sharp\cap B_{r_\sharp}(0)}
v_\eps^\sharp \leq C .
\end{equation}
Set
$$
a_\eps(\zeta):=|\zeta|_{\mathbb H}^{-\eps(Q-2)} .
$$
By~\eqref{eq:Omega-sharp-away-origin}, for every~$\zeta\in\Omega_{\eta_\eps,\rho_\eps}^\sharp$ one has
$$
|\zeta|_{\mathbb H}\ge R_\eps^{-1}.
$$
Therefore
$$
a_\eps(\zeta) \le R_\eps^{\eps(Q-2)}
\leq C\rho_\eps^{-\eps(Q-2)}
\leq C (M_\eps^\eps)^2 .
$$

The bound~\eqref{bound_M_eps} gives
\begin{equation}\label{eq:kelvin-coefficient-uniform}
\|a_\eps\|_{L^\infty(\Omega_{\eta_\eps,\rho_\eps}^\sharp)}
\leq C .
\end{equation}

Next we control the critical mass of $v_\eps^\sharp$ near the origin. By Proposition~\ref{prop:H-Kelv-iso-CR-lap}-(i) and by~\eqref{eq:uniform-tail-v-proof}, for every~$r>0$,
\begin{equation}\label{eq:kelvin-tail-critical}
\int_{\Omega_{\eta_\eps,\rho_\eps}^\sharp\cap B_{2r}(0)}
(v_\eps^\sharp)^{2^\ast}\,{\rm d}\zeta
= \int_{\Omega_{\eta_\eps,\rho_\eps}\setminus B_{\frac{1}{2r}}(0)}
v_\eps^{2^\ast}\,{\rm d}\zeta .
\end{equation}
Fix $p_{\rm o}\in(2,2^\ast)$ and take $\eps_{\rm o} \equiv \eps_{\rm o}(n)>0$ so small that
$$
p_{\rm o}\leq 2^\ast-\eps\leq2^\ast \quad \text{for every }0<\eps<\eps_{\rm o}.
$$
Let $\upsilon_{\rm o}>0$ and $r_{\rm o}>0$ be the constants in Lemma~\ref{lemma:gain}, corresponding to this choice of $p_{\rm o}$, to $q=2^\ast-1$, and to the uniform bound in~\eqref{eq:kelvin-coefficient-uniform}. Set
$$
\vartheta_{\rm o}:=\min\{1,\upsilon_{\rm o}^{2^\ast}\},
$$
and choose $\delta_{\rm o}\in(0,1]$ so small that
$$
2\delta_{\rm o}^{p_{\rm o}/2^\ast}\leq \vartheta_{\rm o}.
$$
By~\eqref{eq:uniform-tail-v-proof} and~\eqref{eq:kelvin-tail-critical}, we can choose
$0<r\leq r_{\rm o}$, independent of~$\eps$, such that, for every sufficiently small~$\eps<\eps_{\rm o}$,
$$
\int_{\Omega_{\eta_\eps,\rho_\eps}^\sharp\cap B_{2r}(0)}
(v_\eps^\sharp)^{2^\ast}\,{\rm d}\zeta
\leq \delta_{\rm o}.
$$
Moreover, decreasing $\eps_{\rm o}$ if necessary, we may also assume that
$$
|\Omega_{\eta_\eps,\rho_\eps}^\sharp\cap B_{2r}(0)|^{\eps/2^\ast}\leq2
\quad \text{for every }0<\eps<\eps_{\rm o}.
$$

Therefore, by H\"older's inequality, 
$$
\int_{\Omega_{\eta_\eps,\rho_\eps}^\sharp\cap B_{2r}(0)} (v_\eps^\sharp)^{2^\ast-\eps}\,{\rm d}\zeta
\leq \vartheta_{\rm o}\,.
$$
Since $\vartheta_{\rm o}\leq \upsilon_{\rm o}^{2^\ast-\eps}$, this gives
\begin{equation}\label{eq:vsharp-critical}
\left( \int_{\Omega_{\eta_\eps,\rho_\eps}^\sharp\cap B_{2r}(0)}
(v_\eps^\sharp)^{2^\ast-\eps}\,{\rm d}\zeta \right)^{\frac{1}{2^\ast-\eps}}
\leq
\upsilon_{\rm o}.
\end{equation}

Applying Lemma~\ref{lemma:gain} to~\eqref{eq:kelvin-v-eps}, and using~\eqref{eq:vsharp-critical} and~\eqref{eq:kelvin-coefficient-uniform}, we obtain
\begin{equation}\label{eq_han4}
\|v_\eps^\sharp\|_{L^{\frac{(2^\ast)^2}{2}}(\Omega_{\eta_\eps,\rho_\eps}^\sharp\cap B_r(0))}
\leq C .
\end{equation}

Now fix $q_0>Q$ such that
\begin{equation}\label{eq:q0-choice-kelvin}
\frac{q_0}{2}(2^\ast-2)
< \frac{(2^\ast)^2}{2}.
\end{equation}
For $\eps$ small enough, we then have
\begin{equation}\label{eq:han-sharp-exp}
\frac{q_0}{2}(2^\ast-\eps-2) \leq \frac{(2^\ast)^2}{2}.
\end{equation}
Define
$$
h_\eps(\zeta) := |\zeta|_{\mathbb H}^{-\eps(Q-2)}
\big(v_\eps^\sharp(\zeta)\big)^{2^\ast-\eps-2}.
$$
By~\eqref{eq:kelvin-coefficient-uniform}, \eqref{eq:han-sharp-exp},~\eqref{eq_han4}, and the fact that~$\Omega_{\eta_\eps,\rho_\eps}^\sharp\cap B_r(0)$ has uniformly bounded measure, we get
\begin{equation}\label{eq:h-eps-q0-bound}
\|h_\eps\|_{L^\frac{q_0}{2}(\Omega_{\eta_\eps,\rho_\eps}^\sharp\cap B_r(0))}
\leq C .
\end{equation}

Writing~\eqref{eq:kelvin-v-eps} as
\[
-\Delta_Xv_\eps^\sharp
=
h_\eps v_\eps^\sharp ,
\]
Lemma~\ref{lemma:sup}, applied in $\Omega_{\eta_\eps,\rho_\eps}^\sharp$ at scale $r/2$ with exponent $q_0>Q$, gives
\[
\sup_{\Omega_{\eta_\eps,\rho_\eps}^\sharp\cap B_\frac{r}{2}(0)}
v_\eps^\sharp
\leq
C
\left(
\frac{1}{|B_r(0)|}
\int_{\Omega_{\eta_\eps,\rho_\eps}^\sharp\cap B_r(0)}
(v_\eps^\sharp)^{2^\ast}\,{\rm d}\zeta
\right)^{1/2^\ast}
\leq C .
\]
Thus~\eqref{boun_max_seq3} holds with $r_\sharp:=r/2$. On the other hand, by~\eqref{eq:kelvin-rough-bound},
\[
v_\eps^\sharp(\zeta)
\leq
|\zeta|_{\mathbb H}^{2-Q}
\leq
r_\sharp^{2-Q}
\quad
\text{for every }
\zeta\in\Omega_{\eta_\eps,\rho_\eps}^\sharp\setminus B_{r_\sharp}(0).
\]
Consequently,
\begin{equation}\label{eq:vsharp-global-bound}
\|v_\eps^\sharp\|_{L^\infty(\Omega_{\eta_\eps,\rho_\eps}^\sharp)}
\leq C .
\end{equation}

We now return to $v_\eps$. For $\zeta\neq0$, the definition of the Kelvin transform gives
$$
v_\eps(\zeta) = |\zeta|_{\mathbb H}^{2-Q}
v_\eps^\sharp(\kappa_{\mathbb H}(\zeta)).
$$
Hence, by~\eqref{eq:vsharp-global-bound},

\begin{equation}\label{perdecay}
v_\eps(\zeta)
\leq C|\zeta|_{\mathbb H}^{2-Q}
\quad \text{for every }\zeta\in\Omega_{\eta_\eps,\rho_\eps}\setminus\{0\}.
\end{equation}

Since~$U$ is positive and continuous, for every fixed~$R>0$ there exists~$c_R>0$ such that
$$
U(\zeta)\geq c_R
\quad\text{for }|\zeta|_{\mathbb H}\leq R.
$$

Moreover, by the explicit Jerison--Lee profile,
$$
U(\zeta)\asymp |\zeta|_{\mathbb H}^{2-Q}
\quad \text{as }|\zeta|_{\mathbb H}\to+\infty .
$$
Combining this lower bound for $U$ with~\eqref{perdecay} in the exterior region, and using $0\leq v_\eps\leq1$ in the bounded region, we obtain
$$
v_\eps(\zeta)\leq C U(\zeta)
\quad \text{for every}~\zeta\in\Omega_{\eta_\eps,\rho_\eps} .
$$
This proves~\eqref{boun_max_seq2} and concludes the proof of Theorem~\ref{han}.
\hfill$\square$

\medskip

\section{Proof of the asymptotic result in Theorem~\ref{cor:centered-koranyi-ball}}\label{sec_localization}

\subsection{Green integral identity}
Before proving Theorem~\ref{cor:centered-koranyi-ball} we establish some integral identity involving the Green and Robin function of a {Kor\'anyi} ball. For this, let us recall the usual surface mean value formula for harmonic function. For this, let~$\Om\subseteq \h$ and~$u \in C^2(\Om)$ be harmonic in~$\Om$, i.~\!e.~$\Delta_X u =0$ in~$\Om$. Then for any~$B_\rho(\xi) \subset \Om$  it holds
\begin{equation}\label{eq:surface-mean}
u(\xi) \ = \  \frac{Q-2}{C_Q \rho^{Q-1} }\int_{\partial B_\rho(\xi)} u(\eta) \frac{|\Dh (|\xi^{-1}\circ \eta|_{\mathbb{H}}|)|^2}{|\nabla(|\xi^{-1}\circ \eta|_{\mathbb{H}}|)|}\,{\rm d}\mathcal{H}^{Q-2}\,,
\end{equation}
where~$C_Q>0$ is the same constant appearing in~\eqref{eq:fundamental}; { see, e.~\!g.,~\cite[Theorem~5.5.4]{BLU08}.}

Moreover, we recall the explicit expression of the Green function with pole at the center of a Kor\'anyi ball~$B_R(\eta)$:
\begin{equation}\label{eq:green-robin-expression}
G_{B_R(\eta)}(\xi;\eta)
:=
K(\eta^{-1}\circ\xi)-\frac{1}{C_QR^{Q-2}}.
\end{equation}
Equivalently, writing
$$
G_{B_R(\eta)}(\xi;\eta)
=
K(\eta^{-1}\circ\xi)-\mathcal{R}_{B_R(\eta)}(0),
$$
with~$\mathcal{R}_{B_R(0)}(0)$ be its regular part at the center.

The following lemma holds true.
\begin{lemma}\label{lemma:green-identity}
For~$R>0$,
\begin{equation}\label{eq:robin-green-2}
\int_{\partial B_R(0)}
|\Dh G_{B_R(0)}(\xi;0)|^2
\langle\mathdutchcal D,n_{\partial B_R(0)}\rangle
\,{\rm d}\mathcal H^{Q-2}
=
(Q-2)\mathcal R_{B_R(0)}(0).
\end{equation}
\end{lemma}

\begin{proof}
For any~$\rho \in (0,R)$ we apply~\cite[Corollary~3.3]{GV00} in~$A_{\rho,R}:= B_R(0)\setminus B_\rho(0)$ to~$G_{B_R(0)}(\xi;0)$ so that
\begin{equation}\label{eq:pohozaev-1}
\begin{aligned}
& 2 \sum_{j=1}^{2n}\int_{\partial A_{\rho,R} }\mathdutchcal{D}G_{B_R(0)}(\xi;0) X_j G_{B_R(0)}(\xi;0) \langle X_j,{n}_{\partial A_{\rho,R}}\rangle \,{\rm d}\mathcal{H}^{Q-2} \\
& \quad +(Q-2)\int_{A_{\rho, R}} \snr{\Dh G_{B_R(0)}(\xi;0)}^2 \,{\rm d}\xi =  \int_{\partial A_{\rho, R}}\snr{\Dh G_{B_R(0)}(\xi;0)}^2 \langle \mathdutchcal{D},{n}_{\partial A_{\rho, R}}\rangle \,{\rm d}\mathcal{H}^{Q-2}\,.
\end{aligned}
\end{equation}

Since~$G_{B_R(0)}(\cdot;0)$ {satisfies
$$
    \Delta_X G_{B_R(0)}(\cdot;0) = 0 \ \ \text{in}~A_{\rho,R} \quad \text{and} \quad 
    G_{B_R(0)}(\cdot;0) = 0 \ \ \text{in}~\partial B_{R}(0)\,,
$$
while~$n_{\partial A_{\rho,R}}= - n_{\partial B_\rho(0)}$ on~$\partial B_\rho(0)$}, integrating by parts in the second integral gives
\begin{multline}\label{eq:pohozaev-2}    
\int_{A_{\rho, R}} \snr{\Dh G_{B_R(0)}(\xi;0)}^2 \,{\rm d}\xi 
\\ =  - \sum_{j=1}^{2n} \int_{\partial B_{\rho}(0)} G_{B_R(0)}(\xi;0) X_j G_{B_R(0)}(\xi;0) \langle X_j,{n}_{\partial B_\rho(0)}\rangle \,{\rm d}\mathcal{H}^{Q-2}.
\end{multline}

{
Moreover, since on~$\partial B_R(0)$ one has~$G_{B_R(0)}(\cdot;0)=0$, the Euclidean gradient of~$G_{B_R(0)}(\cdot;0)$ is normal to~$\partial B_R(0)$ and, for some scalar function~$g$,
$$
\nabla G_{B_R(0)}(\xi;0)=g(\xi)\,{n}_{\partial B_R(0)}
\quad\text{on}~\partial B_R(0).
$$
Consequently,
$$
\sum_{j=1}^{2n}
\mathdutchcal{D}G_{B_R(0)}(\xi;0)\,
X_jG_{B_R(0)}(\xi;0)
\langle X_j,{n}_{\partial B_R(0)}\rangle
=
|\Dh G_{B_R(0)}(\xi;0)|^2
\langle\mathdutchcal{D},{n}_{\partial B_R(0)}\rangle
\quad\text{on } \partial B_R(0).
$$
Combining this identity with~\eqref{eq:pohozaev-2} and~\eqref{eq:pohozaev-1} yields
}
\begin{equation}\label{eq:phozaev-final}
\begin{aligned}
& \int_{\partial B_R(0) }\snr{\Dh G_{B_R(0)}(\xi;0) }^2 \langle \mathdutchcal{D},{n}_{\partial B_R(0)}\rangle \,{\rm d}\mathcal{H}^{Q-2}\\
& \quad =   (Q-2) \sum_{j=1}^{2n} \int_{\partial  B_\rho(0)} \big(X_j G_{B_R(0)}(\xi;0) \big) G_{B_R(0)}(\xi;0) \langle X_j,{n}_{\partial B_\rho(0)}\rangle \,{\rm d}\mathcal{H}^{Q-2}\\
&\qquad - \int_{\partial B_\rho(0)}\snr{\Dh G_{B_R(0)}(\xi;0) }^2 \langle \mathdutchcal{D},{n}_{\partial B_\rho(0)}\rangle \,{\rm d}\mathcal{H}^{Q-2}\\
& \qquad  +\, 2 \sum_{j=1}^{2n}\int_{\partial B_\rho(0) }\mathdutchcal{D}G_{B_R(0)}(\xi;0) X_j G_{B_R(0)}(\xi;0)\langle X_j,{n}_{\partial B_\rho(0)}\rangle \,{\rm d}\mathcal{H}^{Q-2}.
\end{aligned}
\end{equation}

We estimate the right-hand side of~\eqref{eq:phozaev-final} { via~\eqref{eq:surface-mean}}. For this, let us note that by~\eqref{eq:surface-mean} (with~$u=1$) {we can write the first integral as follows}
$$
\begin{aligned}
    & (Q-2) \sum_{j=1}^{2n} \int_{\partial  B_\rho(0)} \big(X_j G_{B_R(0)}(\xi;0) \big) G_{B_R(0)}(\xi;0) \langle X_j,{n}_{\partial B_\rho(0)}\rangle \,{\rm d}\mathcal{H}^{Q-2}\\
    & \quad =  -\frac{1}{\rr^{Q-1}}\left(\frac{Q-2}{C_Q}\right)^2\int_{\partial B_\rho(0)}\left(\frac{1}{\rho^{Q-2}}-\frac{1}{R^{Q-2}}\right)\frac{\snr{\Dh (|\xi|_{\mathbb{H}})}^2}{\snr{\nabla (|\xi|_{\mathbb{H}})}} \,{\rm d}\mathcal{H}^{Q-2} \\
    & \quad = -\frac{Q-2}{C_Q}\left(\frac{1}{\rho^{Q-2}}-\frac{1}{R^{Q-2}}\right).
\end{aligned}
$$
where we have used that by~\eqref{eq:green-robin-expression}  and weak chain-rule
\begin{equation}\label{eq:grad-green}
X_j G_{B_R(0)}(\xi;0) = X_j K(\xi) = -\frac{Q-2}{C_Q \rho^{Q-1}}X_j (|\xi|_{\mathbb{H}})\,,
\end{equation}
and that
\begin{equation}\label{eq:grad-normal} 
    \langle X_j, {n}_{\partial B_\rho (0)}\rangle = \frac{X_j (|\xi|_{\mathbb{H}})}{|\nabla (|\xi|_{\mathbb{H}})|}.
\end{equation}

Moreover, since on~$\partial B_\rho(0)$
$$
|\Dh K(\xi)|^2 = \frac{(Q-2)^2}{C_Q^2\rho^{2Q-2}}|\Dh (|\xi|_{\mathbb{H}})|^2 \quad \text{and} \quad \langle \mathdutchcal{D},{n}_{\partial B_\rho(0)}\rangle = \frac{\rho}{|\nabla (|\xi|_{\mathbb{H}})|}\,,
$$
{given that the gauge~$|\xi|_{\mathbb{H}}$ is $1$-homogeneous,} as before, via~\eqref{eq:surface-mean}  (with~$u=1$), we have that
$$
     \int_{\partial B_\rho(0)}\snr{\Dh G_{B_R(0)}(\xi;0) }^2 \langle \mathdutchcal{D},{n}_{\partial B_\rho(0)}\rangle \,{\rm d}\mathcal{H}^{Q-2} = \frac{Q-2}{C_Q \rho^{Q-2}}.
$$

Since the fundamental solution is homogeneous of degree~$2-Q$, we have~$\mathdutchcal{D} G_{B_R(0)}(\xi;0)= (2-Q)K(\xi)$. Thus, by~\eqref{eq:grad-green},~\eqref{eq:grad-normal} and~\eqref{eq:surface-mean} (with~$u=1$)  we get
$$
\begin{aligned}
     & \sum_{j=1}^{2n}\int_{\partial B_\rho(0) }\mathdutchcal{D}G_{B_R(0)}(\xi;0) X_j G_{B_R(0)}(\xi;0)\langle X_j,{n}_{\partial B_\rho(0)}\rangle \,{\rm d}\mathcal{H}^{Q-2}\\
    & \quad = \frac{1}{\rr^{2Q-3}}\left(\frac{Q-2}{C_Q}\right)^2\int_{\partial B_\rho (0)}\frac{\snr{\Dh (|\xi|_{\mathbb{H}})}^2}{\snr{\nabla |\xi|_{\mathbb{H}}}} \,{\rm d}\mathcal{H}^{Q-2} \\
    & \quad = \frac{Q-2}{C_Q\rho^{Q-2}}
\end{aligned}
$$

Combining all the above information yields~\eqref{eq:robin-green-2}.
\end{proof}

\subsection{Regularity estimates}
In this section we recall and prove some regularity estimates up to the boundary for the horizontal gradient or along the vector fields~$\mathdutchcal{D}$.

\begin{lemma}\label{lemma:reg-ueps}
Let~$u_\eps \in \mathring{\mathcal{D}}^{1,2}_X(B_R(0))$ be a nonnegative weak solution to
    \[
    \begin{cases}
        -\Delta_X u_\eps = u_\eps^{2^\ast-1-\eps} & \text{in}~B_R(0),\\
        u_\eps = 0 & \text{in}~\partial B_R(0).
    \end{cases}
    \]
    Then, there exists~$\eps_0 \equiv \eps_0(Q)>0$ such that, for any~$\eps \in (0,\eps_0)$, the following assertions hold true:
    \begin{enumerate}[\it (i)]
        \item $X_j u_\eps \in L^\infty(B_R(0))$, for every~$1\leq j \leq 2n$;
        \item $\mathdutchcal{D}u_\eps \in L^\infty(B_R(0))$.
        \end{enumerate}
\end{lemma}

\begin{proof}
    The proof follows by the argument developed in~\cite{GV00} and~\cite{Vas06}. See also, Lemma~2.16 and Lemma~2.17 in~\cite{CFL25} for a more recent application in the sub-critical case.
\end{proof}

Before proving the next Lemma we recall the notion of starlikeness with respect to the anisotropic dilations in~\eqref{def_philambda}.

\begin{defn}\label{def:starlike}
Let~$\Omega$ be a $C^1$-domain containing the group identity~$0$. We say that~$\Omega$ is $\delta_\lambda$-starlike with respect to~$0$ along a subset~$K\subseteq\partial\Omega$ if
\[
\langle\mathdutchcal{D},{n}_{\partial\Om}\rangle(\eta)\geq 0
\quad
\text{for every}~\eta\in K.
\]
We say that~$\Omega$ is uniformly $\delta_\lambda$-starlike with respect to~$0$ along~$K$ if there exists~$c_\Omega>0$ such that
\[
\langle\mathdutchcal{D},{n}_{\partial\Om}\rangle(\eta)\geq c_\Omega
\quad
\text{for every}~\eta\in K.
\]
If~$\zeta\in\Omega$, we say that~$\Omega$ is $\delta_\lambda$-starlike, respectively uniformly $\delta_\lambda$-starlike, with respect to~$\zeta$ along~$K$ if the translated domain~$\tau_{\zeta^{-1}}(\Omega)$ has the corresponding property with respect to the identity along~$\tau_{\zeta^{-1}}(K)$.
\end{defn}

\begin{lemma}\label{lemma:reg-ueps-2}
Let~$I_1\Subset I_2$ be neighborhoods of~$\Sigma(B_R(0))=\{(0,\pm R^2)\}$.  Fix~$q_{\rm o}>1$ and~$K_0>0$ and let
$$
q_\eps\in[q_{\rm o},2^\ast-1],
\qquad
a_\eps\in[0,1],
$$
and 
$$
w_\eps\in\mathring{\mathcal D}^{1,2}_X(B_R(0))
$$
be a nonnegative weak solution of
$$
\begin{cases}
-\Delta_Xw_\eps=a_\eps w_\eps^{q_\eps}
& \text{in}~B_R(0),\\
w_\eps=0
& \text{on}~\partial B_R(0).
\end{cases}
$$
Assume that
\begin{equation}\label{eq:w-eps-sup-control}
\sup_\eps
\|w_\eps\|_{L^\infty(B_R(0)\cap I_2)}
\leq K_0.
\end{equation}
Then there exists~$C \equiv C(n,R,I_1,I_2,q_{\rm o},K_0)>0$, 
independent of~$\eps$, such that
\begin{equation}\label{eq:quantitative-characteristic-gradient}
\|\Dh w_\eps\|_{L^\infty(B_R(0)\cap I_1)}
\leq C.
\end{equation}
\end{lemma}

\begin{proof}
Choose fixed neighborhoods~$I_1\Subset J_1\Subset J_2\Subset I_2$ 
of~$\Sigma(B_R(0))$. With the present normalization, the Kor\'anyi ball is the gauge ball considered in~\cite[Theorem~5.3]{GV00} and~\cite[Theorem~5.16]{Vas06}. In particular, it satisfies the $A$-condition (see the definition in Section~5 in~\cite{Vas06}), the convexity condition near its characteristic set, and is uniformly starlike there.

Since~\eqref{eq:w-eps-sup-control} yields
$$
0\leq a_\eps w_\eps^{q_\eps}
\leq
\max\{K_0^{q_{\rm o}},K_0^{2^\ast-1}\}
\quad
\text{in }B_R(0)\cap I_2,
$$
the boundary H\"older estimate under the $A$-condition, used in the proofs of~\cite[Theorem~4.5]{GV00} and~\cite[Theorem~5.6]{Vas06}, gives numbers~$\beta\in(0,1)$ and~$C>0$, independent of~$\eps$, such that
\begin{equation}\label{eq:holder-bd}
w_\eps(\delta_\lambda(\xi))
\leq
C(1-\lambda)^{\beta}\,,
\end{equation}
for every~$\xi$ in a fixed neighborhood of the characteristic set on~$\partial B_R(0)$ and every~$\lambda$ sufficiently close to~$1$. Here, after decreasing~$\beta$ if necessary, we have also used the estimate between the group distance and the dilation parameter proved in~\cite[Lemma~4.1]{GV00}.

We now repeat the barrier iteration in the proof of~\cite[Theorem~4.5]{GV00}. Choose~$N\in\mathbb N$ so large that
$$
\alpha_\ast:=q_{\rm o}^{-N}\leq\beta.
$$
For~$k=0,\dots,N$, set
$$
\alpha_{\eps,k}
:=
\min\{q_\eps^k\alpha_\ast,1\}.
$$
We claim inductively that
\begin{equation}\label{eq:uniform-boundary-decay-iteration}
w_\eps(\delta_\lambda(\xi))
\leq
C_k(1-\lambda)^{\alpha_{\eps,k}},
\end{equation}
where~$C_k$ is independent of~$\eps$. The case~$k=0$ follows from~\eqref{eq:holder-bd}.

Assume~\eqref{eq:uniform-boundary-decay-iteration} holds for some~$k<N$. Since~$q_\eps\in[q_{\rm o},2^\ast-1]$, there exists~$C_k>0$, independent of~$\eps$, such that
$$
a_\eps w_\eps^{q_\eps}
\leq
C_k(1-\lambda)^{q_\eps\alpha_{\eps,k}}
\leq
C_k(1-\lambda)^{\alpha_{\eps,k+1}}.
$$
Let
$$
\Psi_\alpha
:=
(R^4-(|x|^4+t^2))^\alpha
e^{-|x|^2/M},
$$
where~$M>0$ is fixed sufficiently large. The barrier computation in~\cite[Theorem~4.3]{GV00}, equivalently in~\cite[Theorem~5.2]{Vas06}, gives
$$
-\Delta_X\Psi_\alpha
\geq
c\Psi_\alpha
$$
in a fixed dilation collar~$\omega$ of the characteristic set, together with
$$
c_1(1-\lambda)^\alpha
\leq
\Psi_\alpha(\delta_\lambda(\xi))
\leq
c_2(1-\lambda)^\alpha.
$$
We choose~$\omega$ so that~
\begin{enumerate}[\it (i)]
\item $B_R(0)\cap J_2\subset\omega\subset B_R(0)\cap I_2$;
\item every point of~$\omega$ belongs to one of the dilation trajectories used above;
\item $\partial\omega
=
\Gamma_{\partial}\cup\Gamma_{\rm a}$ where~$\Gamma_{\partial}
\subset
\partial B_R(0)$ and~$\Gamma_{\rm a}
\subset \overline{B_R(0)}\cap I_2$ are compact,~$\Gamma_{\partial}$ contains a relative neighborhood of~$\Sigma(B_R(0))$ in~$\partial B_R(0)$, and~$ \Gamma_{\rm a}\cap\partial B_R(0)
\Subset
\partial B_R(0)\setminus\Sigma(B_R(0))$.
\end{enumerate}
Since
$$
\alpha_{\eps,k+1}\in[\alpha_\ast,1],
$$
the constants in the barrier estimates can be chosen uniformly in~$\eps$ and~$k$. Moreover, since every point of~$\omega$ lies on one of the dilation trajectories considered above, there exists~$C_k'>0$, independent of~$\eps$, such that
$$
a_\eps w_\eps^{q_\eps}
\leq
C_k'\Psi_{\alpha_{\eps,k+1}}
\quad
\text{in }\omega.
$$

We next verify the boundary comparison on~$\Gamma_{\rm a}$. Let~$V_{\rm a}\Subset I_2$ be a fixed neighborhood of~$\Gamma_{\rm a}\cap\partial B_R(0)$ whose boundary part is non-characteristic. The uniform bounds for~$w_\eps$ and~$a_\eps w_\eps^{q_\eps}$, together with Theorem~\ref{boundary_hld}, give
$$
\|\Dh w_\eps\|_{L^\infty(B_R(0)\cap V_{\rm a})}
\leq C,
$$
where~$C$ is independent of~$\eps$.

Set
$$
s(x,t)
:=
R^4-(|x|^4+t^2).
$$
Since~$\Dh s\neq0$ on the fixed non-characteristic set
$\Gamma_{\rm a}\cap\partial B_R(0)$, a finite covering by horizontal coordinate neighborhoods and integration along horizontal vector fields transverse to~$\partial B_R(0)$ give
$$
w_\eps(\zeta)
\leq
C s(\zeta)
\quad
\text{for every }
\zeta\in B_R(0)\cap V_{\rm a},
$$
with~$C$ independent of~$\eps$. After decreasing~$V_{\rm a}$, we may assume that~$0\leq s\leq1$ there. Hence, for every~$\alpha\in[\alpha_\ast,1]$,
$$
w_\eps(\zeta)
\leq
Cs(\zeta)
\leq
Cs(\zeta)^\alpha
\leq
C\Psi_\alpha(\zeta)
\quad
\text{in }B_R(0)\cap V_{\rm a}.
$$

On the other hand, on~$\Gamma_{\rm a}\setminus V_{\rm a}
\Subset
B_R(0)\cap I_2$,
$$
\inf_{\substack{
\zeta\in\Gamma_{\rm a}\setminus V_{\rm a}\\
\alpha\in[\alpha_\ast,1]
}}
\Psi_\alpha(\zeta)>0.
$$
Using~\eqref{eq:w-eps-sup-control}, we conclude that there exists~$B_k>0$, independent of~$\eps$, such that
$$
B_k\Psi_{\alpha_{\eps,k+1}}
\geq
w_\eps
\quad
\text{on }\Gamma_{\rm a}.
$$

Choose~$A_k\geq B_k$, independent of~$\eps$, so large that
$$
A_kc\geq C_k'.
$$
On~$\Gamma_{\partial}$ both~$\Psi_{\alpha_{\eps,k+1}}$ and~$w_\eps$ vanish. Therefore
$$
A_k\Psi_{\alpha_{\eps,k+1}}-w_\eps
\geq0
\quad
\text{on }\partial\omega,
$$
and
$$
-\Delta_X
\big(
A_k\Psi_{\alpha_{\eps,k+1}}-w_\eps
\big)
\geq
A_kc\Psi_{\alpha_{\eps,k+1}}
-
a_\eps w_\eps^{q_\eps}
\geq0
\quad
\text{in }\omega.
$$
The weak maximum principle yields
$$
w_\eps
\leq
A_k\Psi_{\alpha_{\eps,k+1}}
\leq
C_{k+1}(1-\lambda)^{\alpha_{\eps,k+1}}
\quad
\text{in }\omega,
$$
where~$C_{k+1}$ is independent of~$\eps$. This proves the induction step. Since
$$
q_\eps^N\alpha_\ast
\geq
q_{\rm o}^N\alpha_\ast
=1,
$$
we obtain
\begin{equation}\label{eq:uniform-linear-boundary-decay}
w_\eps(\xi)
\leq
C\,d_X(\xi)
\quad
\text{in }B_R(0)\cap J_2,
\end{equation}
where
$$
d_X(\xi)
:=
\inf_{\eta\in\partial B_R(0)}
|\eta^{-1}\circ\xi|_{\mathbb H}.
$$
Indeed, the dilation parameter~$1-\lambda$ is uniformly comparable with the Euclidean distance to the boundary near the characteristic set, and the latter is bounded from above by a fixed multiple of~$d_X$.

We now follow the decomposition argument in the proofs of~\cite[Theorem~4.6]{GV00} and~\cite[Theorem~5.10]{Vas06}. Let
$$
\chi\in C^\infty_0(I_2),
\qquad
0\leq\chi\leq1,
\qquad
\chi\equiv1
\quad\text{in }J_2,
$$
and extend
$$
F_\eps
:=
\chi a_\eps w_\eps^{q_\eps}
\mathbbm 1_{B_R(0)}
$$
by zero to~$\h$. Define
$$
v_\eps(\xi)
:=
\int_{\h}
K(\eta^{-1}\circ\xi)F_\eps(\eta)
\,{\rm d}\eta.
$$
Then
$$
-\Delta_Xv_\eps=F_\eps
\quad\text{in }\h.
$$
Since~$\{F_\eps\}$ is uniformly bounded by~\eqref{eq:w-eps-sup-control} and has support in a fixed compact set, the explicit expression of~$K$ and the estimate
$$
|\Dh K(\zeta)|
\leq
C|\zeta|_{\mathbb H}^{1-Q}
$$
give
\begin{equation}\label{eq:uniform-potential-bound}
\|v_\eps\|_{L^\infty(J_2)}
+
\|\Dh v_\eps\|_{L^\infty(J_2)}
\leq C,
\end{equation}
where~$C$ is independent of~$\eps$. 

Set
$$
h_\eps:=w_\eps-v_\eps.
$$
Since~$\chi\equiv1$ in~$J_2$, one has
$$
-\Delta_Xh_\eps=0
\quad
\text{in }B_R(0)\cap J_2.
$$
Let~$\xi\in B_R(0)\cap J_1$ be sufficiently close to~$\partial B_R(0)$. Choose
$$
r=c\,d_X(\xi),
$$
with~$c>0$ sufficiently small that
$$
B_{2r}(\xi)\Subset B_R(0)\cap J_2.
$$
The interior gradient estimate for the harmonic function
$h_\eps-h_\eps(\xi)$ gives
$$
|\Dh h_\eps(\xi)|
\leq
\frac{C}{r}
\sup_{\eta\in B_{2r}(\xi)}
|h_\eps(\eta)-h_\eps(\xi)|.
$$
By~\eqref{eq:uniform-linear-boundary-decay}, for
$\eta\in B_{2r}(\xi)$,
$$
w_\eps(\eta)+w_\eps(\xi)
\leq Cr.
$$
Moreover,~\eqref{eq:uniform-potential-bound} and the intrinsic Poincar\'e estimate give
$$
|v_\eps(\eta)-v_\eps(\xi)|
\leq
C|\xi^{-1}\circ\eta|_{\mathbb H}
\leq Cr.
$$
Consequently,
$$
\sup_{\eta\in B_{2r}(\xi)}
|h_\eps(\eta)-h_\eps(\xi)|
\leq Cr,
$$
and therefore
$$
|\Dh h_\eps(\xi)|\leq C.
$$
Together with~\eqref{eq:uniform-potential-bound}, this yields
$$
|\Dh w_\eps(\xi)|\leq C.
$$
On the remaining compact portion of~$B_R(0)\cap I_1$, which is uniformly separated from~$\partial B_R(0)$, the same estimate follows directly from the interior estimates at a fixed scale. This proves~\eqref{eq:quantitative-characteristic-gradient}.
\end{proof}

\subsection{Pohozaev identity}
We apply the results contained in the previous section to establish a proper Pohozaev identity for~$u_\eps$.

\begin{lemma}\label{lemma:pohozaev-sub}
Let~$u_\eps\in\mathring{\mathcal D}^{1,2}_X(B_R(0))$ be a nonnegative weak solution of
$$
\begin{cases}
-\Delta_Xu_\eps=u_\eps^{2^\ast-\eps-1}
& \text{in }B_R(0),\\
u_\eps=0
& \text{on }\partial B_R(0).
\end{cases}
$$
Assume that~$\Dh u_\eps,\mathdutchcal Du_\eps\in L^\infty(B_R(0))$. Then
\begin{equation}\label{pohozaev_sub_approx}
\frac{\eps(Q-2)}{2^\ast-\eps}
\int_{B_R(0)}u_\eps^{2^\ast-\eps}\,{\rm d}\xi
=
\int_{\partial B_R(0)\setminus\{(0,\pm R^2)\}}
|\Dh u_\eps|^2
\langle\mathdutchcal D,n_{\partial B_R(0)}\rangle
\,{\rm d}\mathcal H^{Q-2}.
\end{equation}
\end{lemma}

\begin{proof}
We argue by exhaustion. For~$r_i\downarrow0$, remove from~$B_R(0)$ two Euclidean neighborhoods of radius~$r_i$ centered at~$(0,\pm R^2)$ and smooth the resulting corners inside annuli of width~$o(r_i)$. This gives smooth domains~$B_i\subset B_R(0)$ with
$$
B_i\uparrow B_R(0),
\qquad
\partial B_i=\gamma_i^{(1)}\cup\gamma_i^{(2)},
$$
where~$\gamma_i^{(1)}\subset\partial B_R(0)\setminus\{(0,\pm R^2)\}$ increases to the non-characteristic boundary and
$$
\mathcal H^{Q-2}(\gamma_i^{(2)})\to0.
$$
Interior hypoelliptic regularity on~$\gamma_i^{(2)}$, non-characteristic boundary regularity on~$\gamma_i^{(1)}$, and a standard regularization justify applying the Pohozaev identity on~$B_i$.

Applying Lemma~\ref{pohozaev} to~$u_\eps$ in~$B_i$, with~$
F(s)={s^{2^\ast-\eps}}/(2^\ast-\eps)$, since~$u_\eps = 0$ on~$\gamma_i^{(1)}$, gives
\begin{eqnarray*}
&& \frac{\eps(Q-2)}{2^\ast-\eps}\int_{B_i}
u_\eps^{2^\ast-\eps}\,{\rm d}\xi\\
&&\quad  =
2\sum_{j=1}^{2n}
\int_{\partial B_i}
\mathdutchcal{D} u_\eps\,X_ju_\eps
\langle X_j,{n}_{\partial B_i}\rangle
\,{\rm d}\mathcal H^{Q-2}
-
\int_{\partial B_i}
|\Dh u_\eps|^2
\langle\mathdutchcal{D},{n}_{\partial B_i}\rangle
\,{\rm d}\mathcal H^{Q-2}
\\
&&\qquad
+
\frac{2}{2^\ast-\eps}
\int_{\gamma_i^{(2)}}
u_\eps^{2^\ast-\eps}
\langle\mathdutchcal{D},{n}_{\partial B_i}\rangle
\,{\rm d}\mathcal H^{Q-2}
+
(Q-2)\sum_{j=1}^{2n}
\int_{\gamma_i^{(2)}}
u_\eps X_ju_\eps
\langle X_j,{n}_{\partial B_i}\rangle
\,{\rm d}\mathcal{H}^{Q-2}.
\end{eqnarray*}

Since~$u_\eps=0$ on~$\partial B_R(0)$, all tangential derivatives vanish on~$\gamma_i^{(1)}$, and the Euclidean gradient is normal to~$\partial B_R(0)$. Thus, for some scalar function~$g_\eps$, we have~$
\nabla u_\eps=g_\eps{n}_{\partial B_R(0)}$ on~$\gamma_i^{(1)}$, which gives
\begin{equation}\label{eq:good-identity}
\sum_{j=1}^{2n}
\mathdutchcal{D} u_\eps\,X_ju_\eps
\langle X_j,{n}_{\partial B_R(0)}\rangle
=
|\Dh u_\eps|^2
\langle\mathdutchcal{D},{n}_{\partial B_R(0)}\rangle
\quad \text{on}~\gamma_i^{(1)}
\end{equation}

It follows that
\begin{eqnarray*}
&& \frac{\eps(Q-2)}{2^\ast-\eps}
\int_{B_i}u_\eps^{2^\ast-\eps}\,{\rm d}\xi \notag\\
&& \quad = \int_{\gamma_i^{(1)}}
|\Dh u_\eps|^2
\langle\mathdutchcal{D},{n}_{\partial B_R(0)}\rangle
\,{\rm d}\mathcal{H}^{Q-2}
+
2\sum_{j=1}^{2n}
\int_{\gamma_i^{(2)}}
\mathdutchcal{D} u_\eps\,X_ju_\eps
\langle X_j,{n}_{\partial B_i}\rangle
\,{\rm d}\mathcal{H}^{Q-2}\\
&& \qquad
-
\int_{\gamma_i^{(2)}}
|\Dh u_\eps|^2
\langle\mathdutchcal{D},{n}_{\partial B_i}\rangle
\,{\rm d}\mathcal H^{Q-2}
+
\frac{2}{2^\ast-\eps}
\int_{\gamma_i^{(2)}}
u_\eps^{2^\ast-\eps}
\langle\mathdutchcal{D},{n}_{\partial B_i}\rangle
\,{\rm d}\mathcal{H}^{Q-2}\\
&& \qquad 
+
(Q-2)\sum_{j=1}^{2n}
\int_{\gamma_i^{(2)}}
u_\eps X_ju_\eps
\langle X_j,{n}_{\partial B_i}\rangle
\,{\rm d}\mathcal{H}^{Q-2}.
\end{eqnarray*}

Since~$u_\eps$,~$\Dh u_\eps$, and~$\mathdutchcal Du_\eps$ are bounded for the fixed~$\eps$, while the Euclidean normals and~$\mathdutchcal D$ are uniformly bounded on the exhausting boundaries, every integral over~$\gamma_i^{(2)}$ tends to zero. On~$\gamma_i^{(1)}$ the integrand is nonnegative because~$B_R(0)$ is~$\delta_\lambda$-starlike; hence monotone convergence applies. This proves~\eqref{pohozaev_sub_approx}.

\end{proof}

\subsection{Proof of Theorem~\ref{cor:centered-koranyi-ball}}

Le us denote with~$\mathcal U(n)$ the unitary group of~$\mathbb C^n$. A function
$u:B_R(0)\to\r$ is called centered-symmetric if
\begin{equation}\label{eq:centered-symmetry}
u(Az,t)=u(z,t) \qquad \text{for every } A\in\mathcal U(n),
\end{equation}
and
\begin{equation}\label{eq:centered-reflection}
u(\overline z,-t)=u(z,t).
\end{equation}

Let~$u_\eps$ be a family of scalar-normalized maximizers solving
\[
\begin{cases}
-\Delta_Xu_\eps=u_\eps^{\frac{Q+2}{Q-2}-\eps}
& \text{in}~B_R(0),\\
u_\eps=0
& \text{on}~\partial B_R(0).
\end{cases}
\]
Assume that~$u_\eps$ is centered-symmetric and satisfies~\eqref{eq:boundary-exlusion}. Since~$0<\delta<R$, compactness of~$\partial B_R(0)$ and continuity of the Kor\'anyi gauge yield~$\delta_E>0$ such that
$$
\{\xi\in B_R(0):
\operatorname{dist}_{\mathbb R^{2n+1}}(\xi,\partial B_R(0))<\delta_E\}
\subset
B_R(0)\setminus B_{R-\delta}(0).
$$
Hence Theorem~\ref{han} applies. We set
$$
M_\eps:=\|u_\eps\|_{L^\infty(B_R(0))},
\qquad
u_\eps(\eta_\eps)=M_\eps,
\qquad
\rho_\eps:=M_\eps^{-\frac{2^\ast-\eps-2}{2}}.
$$
By Theorem~\ref{han}, after passing to a subsequence,
$$
|\Dh u_\eps|^2\,{\rm d}\xi
\tows
(\Sob)^{-\frac{Q-2}{2}}\boldsymbol\delta_{\xi_{\rm o}}
$$
for some~$\xi_{\rm o}\in B_R(0)$. Since each measure is invariant under all transformations in~\eqref{eq:centered-symmetry}--\eqref{eq:centered-reflection}, so is the limit. Hence~$\xi_{\rm o}$ is fixed by all unitary rotations and by~$(z,t)\mapsto(\overline z,-t)$, which implies
\begin{equation}\label{eq:concentration-at-center-proof}
\xi_{\rm o}=0.
\end{equation}

For the sake of readability, we divide the proof into several steps.

\vspace{2mm}

{\it Step 1. The Green profile is selected.}
\vs

Let~$v_\eps$ be given in~\eqref{han_v_eps} and set
\begin{equation}\label{eq:theta-eps}
\theta_\eps
:=
M_\eps^{\frac{(Q-2)\eps}{2}}.
\end{equation}
By~\eqref{bound_M_eps}, up to subsequences, the sequence~$\{\theta_\eps\}$ converges to~$\theta \in [1,C_\infty]$.

We claim that
\begin{equation}\label{eq:W-measure-limit}
M_\eps u_\eps^{2^\ast-\eps-1}\,{\rm d}\xi
\tows
\theta\,\bar c_n\,\boldsymbol{\delta}_0
\quad
\text{on}~\mathcal{M}(\overline{B_R(0)}).
\end{equation}

where
\[
\bar c_n
:=
\int_{\h}U^{2^\ast-1}\,{\rm d}\zeta .
\]
Indeed, let us fix~$\phi\in C(\overline{B_R(0)})$. Using the change of variables~$\xi:=\eta_\eps\circ\delta_{\rho_\eps}(\zeta)$ {and letting~$B_{\eta_\eps,\rho_\eps}:=\delta_{\frac{1}{\rho_\eps}}(\eta_\eps^{-1}\circ B_R(0))$} we get
{
\begin{equation}\label{eq:limit-1}
    \int_{B_R(0)} \phi(\xi)\big(M_\eps u_\eps^{2^\ast-\eps-1}\big)\,{\rm d}\xi
    =
    \theta_\eps
    \int_{B_{\eta_\eps,\rho_\eps}}
    \phi(\eta_\eps\circ \delta_{\rho_\eps}(\zeta))
    v_\eps(\zeta)^{2^\ast-\eps-1}\,{\rm d}\zeta .
\end{equation}
}

Fix~$\alpha\in\big(Q/(Q-2),2^\ast-1\big)$ and choose~$\eps_{\rm o}\equiv\eps_{\rm o}(Q,\alpha)>0$ such that~$\alpha<2^\ast-\eps-1$ for every~$\eps\in(0,\eps_{\rm o})$.

{
After extending~$v_\eps$ by zero outside~$B_{\eta_\eps,\rho_\eps}$,~\eqref{boun_max_seq2}, since~$0\leq v_\eps\leq1$, gives
$$
v_\eps^{2^\ast-\eps-1}
\leq
C U^\alpha 
\quad
\text{in}~\h  \quad \text{with} \quad U^\alpha \in L^1(\h).
$$
}

Moreover, since~$B_{\eta_\eps,\rho_\eps}$ exhausts~$\h$,~\eqref{eq:v-limit-U-unif}, and
$$
\phi(\eta_\eps\circ\delta_{\rho_\eps}(\zeta))
\to
\phi(0),
$$
dominated convergence in~\eqref{eq:limit-1} gives

$$
 \int_{B_R(0)} \phi(\xi)\big(M_\eps u_\eps^{2^\ast-\eps-1}\big)\,{\rm d}\xi
\to 
\theta\,\bar c_n\,\phi(0),
$$
which proves~\eqref{eq:W-measure-limit}.

{
Now, since the function~$w_\eps:=M_\eps u_\eps$ solves
\begin{equation}\label{eq:w-solve}
\begin{cases}
-\Delta_X w_\eps=M_\eps u_\eps^{2^\ast-\eps-1}
& \text{in}~B_R(0),\\
w_\eps=0
& \text{on}~\partial B_R(0),
\end{cases}
\end{equation}
we use the Green representation formula. Let~$G_{B_R(0)}(\xi;\eta)$ denote the Dirichlet Green kernel of~$-\Delta_X$ in~$B_R(0)$. We shall use the standard comparison estimate
$$
0\leq G_{B_R(0)}(\xi;\eta)
\leq K(\eta^{-1}\circ\xi),
\qquad
\xi,\eta\in B_R(0),\quad \xi\neq\eta.
$$
Then
$$
w_\eps(\xi)
=
\int_{B_R(0)}
G_{B_R(0)}(\xi;\eta)
\big(M_\eps u_\eps^{2^\ast-\eps-1}\big)
\,{\rm d}\eta.
$$

Let~$\mathcal K\Subset\overline{B_R(0)}\setminus\{0\}$ be compact, with~$\mathcal K\cap\partial B_R(0)$ contained in a non-characteristic portion of the boundary. Choose~$r \in (0,R/3)$ such that
$$
\mathcal K\cap\overline{B_{3r}(0)}=\emptyset,
$$
and let~$\chi\in C^\infty_0(B_{2r}(0))$ satisfy~$0\leq\chi\leq1$ and~$\chi\equiv1$ on~$B_r(0)$. Since
$$
\left\{\eta\mapsto
\chi(\eta)G_{B_R(0)}(\xi;\eta):\xi\in\mathcal K\right\}
$$
is compact in~$C(\overline{B_R(0)})$, the weak convergence~\eqref{eq:W-measure-limit} gives
$$
\sup_{\xi\in\mathcal K}
\left|
\int_{B_R(0)}
\chi(\eta)G_{B_R(0)}(\xi;\eta)
\big(M_\eps u_\eps^{2^\ast-\eps-1}\big)
\,{\rm d}\eta
-
\theta\bar c_nG_{B_R(0)}(\xi;0)
\right|
\to0.
$$

On~$B_R(0)\setminus B_r(0)$, since~$\eta_\eps\to0$, for~$\eps$ sufficiently small one has
$$
|\eta_\eps^{-1}\circ\eta|_{\mathbb H}\geq c_r>0.
$$
Thus~\eqref{bound_max_seq},~\eqref{bound_M_eps}, and
$$
\rho_\eps^{Q-2}
=
M_\eps^{-2+\frac{(Q-2)\eps}{2}}
$$
give
\begin{equation}\label{eq:glob-est}
u_\eps(\eta)
\leq
\frac{C_r}{M_\eps}.
\end{equation}
Therefore
\begin{equation}\label{eq:blob-est}
M_\eps u_\eps(\eta)^{2^\ast-\eps-1}
\leq
C_rM_\eps^{-(2^\ast-\eps-2)}
=
\frac{C_rM_\eps^\eps}{M_\eps^{2^\ast-2}}
\leq
\frac{C_r}{M_\eps^{2^\ast-2}}
\to0.
\end{equation}
Using the comparison estimate for the Green kernel and the local integrability of~$K$, we obtain
$$
\sup_{\xi\in\mathcal K}
\int_{B_R(0)}
(1-\chi(\eta))G_{B_R(0)}(\xi;\eta)
\big(M_\eps u_\eps^{2^\ast-\eps-1}\big)
\,{\rm d}\eta
\to0.
$$
Consequently,
\begin{equation}\label{eq:green-profile-theta}
M_\eps u_\eps
\to
\theta\bar c_nG_{B_R(0)}(\cdot;0)
\quad
\text{locally uniformly in }
\overline{B_R(0)}\setminus\{0\}
\end{equation}
on compact sets whose boundary part is non-characteristic.

Let~$\gamma\Subset\partial B_R(0)\setminus\{(0,\pm R^2)\}$ and choose a boundary neighborhood~$V$ of~$\gamma$ with~$0\notin\overline V$. Set
$$
z_\eps
:=
w_\eps-\theta\bar c_nG_{B_R(0)}(\cdot;0).
$$
Then~$z_\eps=0$ on~$\partial B_R(0)\cap V$,
$$
\|z_\eps\|_{L^\infty(B_R(0)\cap V)}\to0,
$$
and
$$
-\Delta_Xz_\eps
=
M_\eps u_\eps^{2^\ast-\eps-1}
\to0
\quad
\text{in }L^\infty(B_R(0)\cap V).
$$
Theorem~\ref{boundary_hld}, applied on a smaller boundary neighborhood, gives
\begin{equation}\label{eq:green-profile-boundary}
M_\eps u_\eps
\to
\theta\bar c_nG_{B_R(0)}(\cdot;0)
\quad
\text{in }\Gamma^{1,\beta}_{\rm loc}
(\partial B_R(0)\setminus\{(0,\pm R^2)\}).
\end{equation}
It remains to control the blow-up near the characteristic set. Let~$I_1\Subset I_\sigma$ be fixed neighborhoods of~$\{(0,\pm R^2)\}$ such that~$0\notin\overline{I_\sigma}$. By~\eqref{eq:glob-est},
$$
\|w_\eps\|_{L^\infty(B_R(0)\cap I_\sigma)}
\leq C.
$$
Set
$$
q_\eps:=2^\ast-\eps-1,
\qquad
a_\eps
:=
M_\eps^{-(2^\ast-\eps-2)}
=
\rho_\eps^2.
$$
Since~$w_\eps=M_\eps u_\eps$, equation~\eqref{eq:w-solve} can be rewritten as
$$
\begin{cases}
-\Delta_Xw_\eps=a_\eps w_\eps^{q_\eps}
& \text{in }B_R(0),\\
w_\eps=0
& \text{on }\partial B_R(0).
\end{cases}
$$
By~\eqref{han_supremum_limit}, for~$\eps$ sufficiently small,
$$
0<a_\eps\leq1.
$$
Moreover, after decreasing~$\eps_{\rm o}$ if necessary, there exists~$q_{\rm o}>1$ such that
$$
q_{\rm o}\leq q_\eps\leq2^\ast-1
\qquad
\text{for every }0<\eps<\eps_{\rm o}.
$$
Applying Lemma~\ref{lemma:reg-ueps-2} on~$I_1\Subset I_\sigma$, we obtain
\begin{equation}\label{eq:Y-boundary-char}
\|\Dh w_\eps\|_{L^\infty(B_R(0)\cap I_1)}
\leq C_\sigma,
\end{equation}
where~$C_\sigma$ is independent of~$\eps$.

\vspace{2mm}
{\it Step 2. Passing to the limit in the Pohozaev identity.}
\vspace{1mm}

Applying Lemma~\ref{lemma:pohozaev-sub}, we multiply~\eqref{pohozaev_sub_approx} by~$M_\eps^2$ obtaining
\begin{equation}\label{eq:pohozaev-multiplied}
\frac{\eps(Q-2)}{2^\ast-\eps}
M_\eps^2
\int_{B_R(0)}u_\eps^{2^\ast-\eps}\,{\rm d}\xi
=
\int_{\partial B_R(0) {\setminus \{(0,\pm R^2)\}}}
|\Dh(M_\eps u_\eps)|^2
\langle\mathdutchcal{D}, {n}_{\partial B_R(0)}\rangle
\,{\rm d}\mathcal H^{Q-2}.
\end{equation}

{
We pass to the limit in the boundary term. Write again~$w_\eps:=M_\eps u_\eps$, and let~$I_\delta\subset I_1$ be a decreasing family of relatively open neighborhoods of~$\{(0,\pm R^2)\}$ such that
$$
\mathcal H^{Q-2}(\partial B_R(0)\cap I_\delta)\to0
\qquad\text{as }\delta\to0^+.
$$
For every fixed~$\delta>0$,~\eqref{eq:green-profile-boundary} gives
$$
\Dh w_\eps
\to
\theta\bar c_n\Dh G_{B_R(0)}(\cdot;0)
$$
uniformly on~$\partial B_R(0)\setminus I_\delta$. Hence
$$
\begin{aligned}
&\lim_{\eps\to0^+}
\int_{\partial B_R(0)\setminus I_\delta}
|\Dh w_\eps|^2
\langle\mathdutchcal D,n_{\partial B_R(0)}\rangle
\,{\rm d}\mathcal H^{Q-2}\\
&\qquad=
\theta^2\bar c_n^2
\int_{\partial B_R(0)\setminus I_\delta}
|\Dh G_{B_R(0)}(\xi;0)|^2
\langle\mathdutchcal D,n_{\partial B_R(0)}\rangle
\,{\rm d}\mathcal H^{Q-2}(\xi).
\end{aligned}
$$
On the other hand,~\eqref{eq:Y-boundary-char} and the boundedness of~$\langle\mathdutchcal D,n_{\partial B_R(0)}\rangle$ give
$$
\sup_{\eps}
\int_{\partial B_R(0)\cap I_\delta}
|\Dh w_\eps|^2
|\langle\mathdutchcal D,n_{\partial B_R(0)}\rangle|
\,{\rm d}\mathcal H^{Q-2}
\leq
C_\sigma^2
\mathcal H^{Q-2}(\partial B_R(0)\cap I_\delta)
\to0
$$
as~$\delta\to0^+$. The same estimate applies to the fixed function~$G_{B_R(0)}(\cdot;0)$. Therefore, first letting~$\eps\to0^+$ with~$\delta$ fixed and then letting~$\delta\to0^+$, we obtain
\begin{equation}\label{eq:phozaev-right}
\begin{aligned}
&\lim_{\eps\to0^+}
\int_{\partial B_R(0)\setminus\{(0,\pm R^2)\}}
|\Dh w_\eps|^2
\langle\mathdutchcal D,n_{\partial B_R(0)}\rangle
\,{\rm d}\mathcal H^{Q-2}\\
&\qquad=
\theta^2\bar c_n^2
\int_{\partial B_R(0)}
|\Dh G_{B_R(0)}(\xi;0)|^2
\langle\mathdutchcal D,n_{\partial B_R(0)}\rangle
\,{\rm d}\mathcal H^{Q-2}.
\end{aligned}
\end{equation}
}

Combining~\eqref{eq:phozaev-right} in~\eqref{eq:pohozaev-multiplied} and letting~$\eps\to0^+$ yields
{
\begin{equation}\label{eq:rate-with-theta}
\lim_{\eps\to0}
\eps M_\eps^2
=
\frac{2^\ast\theta^2\bar c_n^2}{(Q-2)(\Sob)^{-\frac{Q-2}{2}}}
\int_{\partial B_R(0)}
|\Dh G_{B_R(0)}(\xi;0)|^2
\langle\mathdutchcal D,n_{\partial B_R(0)}\rangle
\,{\rm d}\mathcal H^{Q-2}.
\end{equation}
}

{
In particular,~\eqref{eq:rate-with-theta} yields that~$\eps M_\eps^2$ is bounded when~$\eps\to0^+$ along subsequences. Hence, there exists~$C=C(B_R(0),Q,\theta)>0$ such that, as~$\eps\to0^+$,
}
$$
\eps\log M_\eps
\ \leq \ 
\frac{\eps}{2}\log\frac{C}{\eps} \to 0\,,
$$
i.~\!e., recalling~\eqref{eq:theta-eps} we have
$$
\theta_\eps
=
e^{
\frac{(Q-2)\eps}{2}\log M_\eps
}  \to 1.
$$
Thus~$\theta=1$.  Consequently,~\eqref{eq:green-profile-theta} and~\eqref{eq:rate-with-theta} become
{
\begin{equation}\label{eq:green-profile-final-proof}
\begin{aligned}
& M_\eps u_\eps
\to
\bar c_nG_{B_R(0)}(\cdot;0)
\quad\text{locally uniformly in }B_R(0)\setminus\{0\},\\
& \lim_{\eps\to0^+}
\eps M_\eps^2
=
\frac{2^\ast\bar c_n^2}{(Q-2)(\Sob)^{-\frac{Q-2}{2}}}
\int_{\partial B_R(0)}
|\Dh G_{B_R(0)}(\xi;0)|^2
\langle\mathdutchcal D,n_{\partial B_R(0)}\rangle
\,{\rm d}\mathcal H^{Q-2}.
\end{aligned}
\end{equation}
}

{
Therefore, combining~\eqref{eq:green-profile-final-proof} with Lemma~\ref{lemma:green-identity} we obtain
$$
\lim_{\eps\to0^+}\eps M_\eps^2
=
\frac{2^\ast\bar c_n^2}{(\Sob)^{-\frac{Q-2}{2}}}
\frac{1}{C_QR^{Q-2}}.
$$
Moreover, for every~$\xi\in B_R(0)\setminus\{0\}$,
$$
\frac{u_\eps(\xi)}{\sqrt\eps}
=
\frac{M_\eps u_\eps(\xi)}{\sqrt{\eps M_\eps^2}}
\to
\sqrt{\frac{(\Sob)^{-\frac{Q-2}{2}}R^{Q-2}}{2^\ast C_Q}}
\left(
\frac{1}{|\xi|_{\mathbb H}^{Q-2}}
-
\frac{1}{R^{Q-2}}
\right).
$$
Since the preceding argument applies to every sequence~$\eps_j\to0^+$ and every subsequence has a further subsequence with the same limits, the limits above hold for the whole family as~$\eps\to0^+$.
This concludes the proof.
\hfill $\square$
}

\vspace{3mm}

 \mbox{}
 \\ {\bf Conflict of interest.} \, The authors declare to have no conflict of interests. No data are attached to this paper.


\medskip
\begin{thebibliography}{99}
	
	
	\bibitem{AG21} {N. Alamri, N. Gamara}: Non-characteristic Heisenberg group domains. {\it Period. Math. Hung.} {\bf 82} (2021), 16--28.
	\vspace{-2.5mm}

	
	\bibitem{AG03} {M. Amar, A. Garroni}: {$\Gamma$-convergence of concentration problems}. {\it Ann. Scuola Norm. Sup. Pisa Cl. Sci.} {\bf 2}~(2003), 151--179.
	\vs
	
	\bibitem{AP87} {F. V. Atkinson, L. A. Peletier}: Elliptic equations with nearly critical growth. {\it J. Differential Equations} {\bf 70} (1987), no. 3, 349--365.
	\vs
	
	\bibitem{BC88} {A.~Bahri, J.~\!-M. Coron}: On a nonlinear elliptic equation involving the critical Sobolev exponent: the effect of the topology of the domain. {\it Comm. Pure Appl. Math.} {\bf 41} (1988), no. 3, 253--294.
	\vs
	
	
	
	
	\bibitem{BGM19} {A. Banerjee, N. Garofalo, I.~\!H.~Munive}: Compactness methods for $\Gamma^{1,\alpha}$ boundary Schauder estimates in Carnot groups. {\it Calc. Var. Partial Differential Equations} {\bf 58} (2019), no. 3, 29--97.
	\vs
	
	\bibitem{BGM22} {A. Banerjee, N. Garofalo, I.~\!H.~Munive}:
	Higher order Boundary Schauder Estimates in Carnot Groups. {\it Math. Ann.}, to appear.
	\vs
	
	
	
	
	
	
    
	\bibitem{BLU08} {A. Bonfiglioli, E. Lanconelli, F. Uguzzoni}:
	{\it Stratified Lie Groups and their sub-Laplacians}.
	Springer Monographs in Mathematics, Springer, Berlin, 2007. 
	\vs
	
	

	
	\bibitem{BP89}{H. Brezis, L.~\!A. Peletier}: Asymptotic for Elliptic Equations involving critical growth. In {\it Partial Differential Equations and the Calculus of Variations. Essays in Honor of Ennio De Giorgi, Vol. 1}, Progr. Differ. Equ. Appl., Birkh\"auser, Boston (1989), 149--192.
	\vs
	
	
	
	\bibitem{CFL25}{H.~Chen, Y.~\!L.~Fan, X.~Liao} Multiple sign-changing solutions of Brezis-Nirenberg problem for Kohn Laplacian on a partially symmetric domain in Heisenberg group. {\it Calc. Var. Partial Differential Equations}  {\bf 64}, 171 (2025).
	\vspace{-2mm}

	
	\bibitem{CGN02} {L. Capogna, N. Garofalo, D.-M.~Nhieu}: Properties of Harmonic Measures in the Dirichlet Problem for Nilpotent Lie Groups of Heisenberg Type. {\it Amer. J. Math.}~{\bf 124} (2002), no.~2, 273--306.
	\vs


    \bibitem{CGL93} {G. Citti, N. Garofalo, E. Lanconelli}:  Harnack's inequality for sum of squares of vector fields plus a potential. {\it Amer. J. Math.} (3) {\bf 115}~(1993), 69--734.
\vs
	
	
	
	\bibitem{CGS21}{G. Citti, G. Giovannardi, Y. Sire}: Schauder estimates up to the boundary on H-type groups: an approach via the double layer potential. {\it  Ann. Scuola Norm. Sup. Pisa Cl. Sci.}~{\bf 37} (2024), \href{https://doi.org/10.2422/2036-2145.202302_015}{\tt DOI:10.2422/2036-2145.202302\_015}
	\vs
 
	\bibitem{CU01}{G. Citti, F. Uguzzoni}: Critical semilinear equations on the Heisenberg group: the effect of the topology of the domain. {\it Nonlinear Anal.}~{\bf 46}~(2001), 399--417.
	\vs
	
	
	
	
	
	
	\bibitem{DPMP10} {M.~del Pino, M.~Musso, F.~Pacard}:
	Bubbling along boundary geodesics near the second critical exponent. {\it J. Eur. Math. Soc. (JEMS)} {\bf 12}~(2010), No.~6, 1553--1605.
	\vs
	
	
	
	
	
	
	
	
	\bibitem{FGM02}{M. Flucher, A. Garroni, S. M\"uller}: Concentration of low energy extremals: Identification of concentration points. {\it Calc. Var. Partial Differential Equations}~{\bf 14} (2002), 483--516.
	\vs
	
	\bibitem{FM99}{M. Flucher, S. M\"uller}: Concentration of low energy extremals. {\it Ann. Inst. H. Poincar\'e - Anal. Non Lin\'eaire}~{\bf 16} (1999), no.~3, 269--298.
	\vs
	
	
	\bibitem{Fol75}{G.~\!B.~Folland}: Subelliptic estimates and function spaces on nilpotent Lie groups. {\it Ark. Math.} {\bf 13} (1975), 161--207.
	\vs
	

	
	\bibitem{FS74}{G.~\!B.~Folland, E.~\!M. Stein}: Estimates for the $\bar\partial_b$ complex and analysis on the Heisenberg group. {\it Comm. Pure Appl. Math.}~{\bf 27} (1974), 429--522.
	\vs

{  	\bibitem{FS82}{G.~\!B.~Folland, E.~\!M. Stein}: {\it Hardy spaces on homogeneous groups}. Volume 28 of Mathematical Notes. Princeton University Press, Princeton, N.J 1982.}
    \vs
	
	\bibitem{FL12}{R. Frank, E.~\!H. Lieb}: Sharp constants in several inequalities on the Heisenberg group. {\it Ann. Math.}~{\bf 176}~(2012), 349--381. 
	\vs
	
	
	\bibitem{FGMT15} {R.~\!L. Frank, M.~\!d.~\!M. Gonz\'alez, D.~\!D. Monticelli, J. Tan}: An extension problem for the CR fractional Laplacian. {\it Adv. Math.}~{\bf 270} (2015), 97--137.
	\vs
	
	
	\bibitem{FKK23} {R. L. Frank, T. K\"onig and H. Kovarik}: Blow-up of solutions of critical elliptic equation in three dimensions. {\it Anal. PDE}, to appear. 
	\vs
	
	
	
	\bibitem{Gam01}{N. Gamara}: The CR Yamabe conjecture -- the case $n=1$. {\it J. Eur. Math. Soc.} {\bf 3} (2001), no.~2, 105--137.
	\vs
	
	\bibitem{Gam17}{N. Gamara, A. Makni}:
The First Eigenvalue of the Kohn--Laplace Operator in the Heisenberg Group.
{\it Mediterr. J. Math.} {\bf 14} (2017), Art.~60.
\vs
	

	\bibitem{GL92}{N. Garofalo, E. Lanconelli}: Existence and Nonexistence Results for Semilinear Equations on the Heisenberg Group. {\it Indiana Univ. Math. J.} {\bf 41}~(1992), no.~1, 71--98.
	
	\vs
	
	
	\bibitem{GLV23}{N. Garofalo, A. Loiudice, D. Vassylev}:  Optimal decay for solutions of nonlocal semilinear equations with critical exponent in homogeneous group. {\it Proc. Royal Soc. Edinburgh Section~A}, to appear.
	\vs
	
	
	
	
	\bibitem{GV00}{N. Garofalo, D. Vassilev}: Regularity near the characteristic set
	in the non-linear Dirichlet problem
	and conformal geometry of sub-Laplacians on Carnot groups. {\it Math. Ann.} {\bf 318}~(2000), 453--516.
	\vs
	

	
	\bibitem{Han91} {Z.-C. Han}: {Asymptotic approach to singular solutions for nonlinear elliptic equations involving critical Sobolev exponent}. {\it Ann. Inst. Henri Poincar\'e Anal. Non Lin\'eaire} {\bf 8}~(1991),  159--174.
	\vs	
	
	\bibitem{IV11}{S.~\!P. Ivanov, D.~\!B. Vassilev}: {\it Extremals for the Sobolev Inequality and the Quaternionic Contact Yamabe Problem}. World Scientific Publishing, Singapore, 2011.
	
	\vs
	
	\bibitem{Jer81} {D. Jerison}: The Dirichlet problem for the Kohn Laplacian on the Heisenberg group. I. {\it J. Functional Analysis} {\bf43} (1981), no. 1, 97--142.
	\vs
	
	\bibitem{Jer81b} {D. Jerison}: The Dirichlet problem for the Kohn Laplacian on the Heisenberg group. II. {\it J. Functional Analysis} {\bf43}  (1981), no. 2, 224--257.
	\vs
	
	
	\bibitem{JL87} {D. Jerison, J.~\!M. Lee}: The Yamabe problem on CR manifolds. {\it J. Differential Geom.} {\bf 25} (1987), no. 2, 167--197.
		\vspace{-2mm}
 
	
	\bibitem{JL88} {D. Jerison, J.~\!M. Lee}: Extremals for the Sobolev inequality on the Heisenberg group and the~CR~Yamabe problem. {\it J. Amer. Math. Soc.}~{\bf 1}~(1988), no.~1, 1--13.
	\vs
		
	\bibitem{LU98} {E. Lanconelli, F. Uguzzoni}:
	Asymptotic behavior and non-existence theorems for semilinear Dirichlet problems involving critical exponent on unbounded domains of the Heisenberg group.
	{\it Boll. UMI (Serie 8)} {\bf 1-B}~(1998), no.~1, 139--168.
	\vs
	
	
	
	\bibitem{Loi05} {A. Loiudice}: Improved Sobolev inequalities on the Heisenberg group. {\it Nonlinear Analysis: Theory, Methods \& Applications}, {\bf 62}~(2005), no.~5, 953--962.
	\vs
	
	\bibitem{MMP13} {A. Maalaoui, V. Martino, A. Pistoia}: Concentrating solutions for a sub-critical sub-elliptic problem. {\it Diff. Int. Eq.}~{\bf 26}~(2013), no. 11-12, 1263--1274.
	\vs
	
	
	
	\bibitem{MPPP23} {M. Manfredini, G. Palatucci, M. Piccinini, S. Polidoro}:
	{H\"older continuity and boundedness estimates for nonlinear fractional equations in the    Heisenberg group}. {\it J. Geom. Anal.}~{\bf 33}~(2023), no.~3, Art.~77. 
	\vs

	\bibitem{Pal11}{G. Palatucci}: Subcritical approximation of the Sobolev quotient and a related concentration result. {\it Rend. Sem. Mat. Univ. Padova} {\bf 125}~(2011), 1--14.
	\vs
	
	\bibitem{Pal11b}{G. Palatucci}: $p$-Laplacian problems with critical Sobolev exponent.
	{\it Asymptot. Anal.} {\bf 73}~(2011), no.~1-2, 37--52.
	\vspace{-2mm}
	
	\bibitem{PP22} {G. Palatucci, M. Piccinini}:
	{Nonlocal Harnack inequalities in the Heisenberg group}.
	{\it Calc. Var. Partial Differential Equations}~{\bf 61} (2022), Art.~185.
	\vs
	
	\bibitem{PPT25} {G. Palatucci, M. Piccinini, L. Temperini}:
	Struwe's Global Compactness and energy approximation of the critical Sobolev embedding in the Heisenberg group. {\it Adv. Calc. Var.}~{\bf 18} (2025), no.~3, 731--754. 
	\vs
	
	\bibitem{PP23} {G. Palatucci, M. Piccinini, L. Temperini}:
	{Global Compactness, subcritical approximation of the Sobolev quotient, and a related concentration result in the Heisenberg group}. 
	{\it Trends in Mathematics},  Birkh\"auser,   	2024.
	\vs
	
	
	\bibitem{Pas93}{D.  Passaseo}: Nonexistence results for elliptic problems with supercritical nonlinearity in nontrivial domains. {\it J. Funct. Anal.}  {\bf 114} (1993), 97--105.
	\vs
	
	
	
	\bibitem{PR03} {A. Pistoia,  O. Rey}: {Boundary blow-up
		for a Brezis-Peletier problem on a singular domain}.  {\it Calc. Var. Partial Differential Equations}~{\bf 18} (2003), no.~3, 243--251.
	\vs
	
	\bibitem{PT23} {P. Pucci, L. Temperini}: {On the Concentration-compactness principle for Folland-Stein spaces and for fractional horizontal Sobolev spaces}. {\it Math. Eng.}~{\bf 5} (2023), no.~1, 1--21. 
	\vs
	
	\bibitem{Rey89} {O. Rey}: {Proof of the conjecture of H.~Brezis and L.~\!A.~Peletier}. {\it Manuscripta math.} {\bf 65} (1989),  19--37. 
	\vs
	
	
	
	\bibitem{Ugu99} {F. Uguzzoni}:  A non-existence theorem for a semilinear Dirichlet problem involving critical exponent on halfspaces of the Heisenberg group. {\it NoDEA Nonlinear Differential Equations Appl.} {\bf 6} (1999), no. 2, 191--206.
		\vspace{-2mm}
 
	
	\bibitem{Vas06} {D. Vassilev}: Existence of solutions and regularity near the characteristic boundary for sub-Laplacian equations on Carnot groups. {\it Pacific J. Math.} {\bf 227} (2006), no. 2, 361--397. 
	\vs
	
	\bibitem{Wei98}{J. Wei}: Asymptotic behavior of least-energy solutions of a semilinear Dirichlet problem involving critical Sobolev exponent. {\it J. Math. Soc. Japan}~{\bf 1}~(1998), 139--153.
	\vs
	
	
	
	
	
\end{thebibliography}
\end{document}